\documentclass[12pt,notitlepage,twoside]{amsart}
\usepackage[T1]{fontenc}
\usepackage{a4}
\usepackage{amscd}
\usepackage[dvips]{color,graphicx}
\usepackage{epsfig}
\numberwithin{equation}{section}

\author{Jörn Peter}
\title{Hausdorff measure of Julia sets in the exponential family}
\address{Mathematisches Seminar, Christian-Albrechts-Universität zu Kiel, 24098 Kiel, Germany}
\curraddr{Departament de Matemàtica Aplicada i Anàlisi, Universitat de Barcelona, 08007 Barcelona, Spain}
\email{peter@math.uni-kiel.de}

\newcommand{\abs}[1]{\left|#1\right|}

\newcommand{\trenn}{\hspace{0.1cm}|\hspace{0.1cm}}
\newcommand{\Trenn}{\hspace{0.1cm}{\Big|}\hspace{0.1cm}}
\newcommand{\set}[1]{\left\{#1\right\}}
\newcommand{\bra}[1]{\left(#1\right)}
\newcommand{\Com}{\mathbb C}

\newcommand{\eps}{\varepsilon}

\newtheoremstyle{myintro}
{10pt}{10pt}{\itshape}{}{\bfseries}{}{ }{\thmname{#1}\thmnote{
#3}}

\newtheoremstyle{myplain}
{10pt}{10pt}{\itshape}{}{\bfseries}{}{ }{\thmname{#1}\thmnumber{
#2}\thmnote{ (#3)}}

\newtheoremstyle{mydefinition}
{10pt}{10pt}{}{}{\bfseries}{}{ }{\thmname{#1}\thmnumber{
#2}\thmnote{ (#3)}}

\theoremstyle{myplain}
\newtheorem{theor}{Theorem}[section]
\newtheorem{lemmm}[theor]{Lemma}
\newtheorem{coror}[theor]{Corollary}

\theoremstyle{mydefinition}
\newtheorem{defin}[theor]{Definition}

\theoremstyle{myintro}
\newtheorem{theorin}[theor]{Theorem}

\begin{document}

\newenvironment{remark}{\noindent \textbf{Remark:}}{\hfill$\spadesuit$}
\newenvironment{remarks}{\noindent \textbf{Remarks:}}{\hfill$\spadesuit$}
\newenvironment{proo}{\noindent \textbf{Proof.}}{\hfill$\Box$}
\newenvironment{proo2}{\noindent \textbf{Proof of Theorem \ref{zeromeas}.}}{\hfill$\Box$}
\newenvironment{proo3}{\noindent \textbf{Proof of Theorem \ref{main1}.}}{\hfill$\Box$}
\newenvironment{theo}{\begin{theor}}{\end{theor}}
\newenvironment{lemm}{\begin{lemmm}}{\end{lemmm}}
\newenvironment{coro}{\begin{coror}}{\end{coror}}
\newenvironment{defi}{\begin{defin}}{\end{defin}}
\newenvironment{theoin}{\begin{theorin}}{\end{theorin}}

\maketitle

\begin{abstract}
We consider the Hausdorff measure of Julia sets and escaping sets of exponential maps with respect to certain
gauge functions. We give conditions on the growth of the gauge function which imply that the measure is
zero or infinity, respectively.
\end{abstract}

\section{Main result and outline}

\subsection{Introduction and main result}

The exponential family consists of all
functions $E_\lambda(z):=\lambda e^z$, where $\lambda\in\Com\setminus\set{0}$. Denote by $E_\lambda^n$ the $n$-th iterate of $E_\lambda$.
The \emph{Fatou set} $\mathcal{F}(E_\lambda)$ is the set of points where the iterates
$E_\lambda^n$ of $E_\lambda$ form a normal family in the sense of Montel (or equivalently, where the iterates
are equicontinuous). The complement of $\mathcal{F}(E_\lambda)$ is the \emph{Julia set}
$\mathcal{J}(E_\lambda)$. The \emph{escaping set} $I(E_\lambda)$ is the set of all points $z$ such that $E_\lambda^n(z)$
tends to infinity as $n\to\infty$. A result
of Eremenko and Lyubich \cite{eremenkolyubich} implies that $I(E_\lambda)\subset\mathcal{J}(E_\lambda)$.
The function $E_\lambda$ is called \emph{hyperbolic} if there exists $z_0\in\Com$ with $E_\lambda^n(z_0)=z_0$ and
$\abs{(E_\lambda^n)'(z_0)}<1$.

A \emph{gauge function} is a monotonically
increasing function $h:[0,\eps)\to\mathbb{R}_{\geq 0}$ (where $\eps>0$) which is
continuous from the right and satisfies $h(0)=0$. Define
$$\mathcal{H}^h(A):=\lim_{\delta\to 0}\inf\set{\sum_{i=1}^\infty h(\text{diam
}A_i)\Trenn\bigcup_{i=1}^\infty A_i\supset A,\text{ diam
}A_i<\delta\text{ for every }i}.$$
Then $\mathcal{H}^h$ is a metric outer measure on $\Com$, called
the \emph{Hausdorff measure with respect to $h$}. In the special case where $h^s(t):=t^s$ for some $s>0$, $\mathcal{H}^{h^s}$ is the
\emph{$s$-dimensional outer Hausdorff measure}. Given $A\subset\Com$, it is well known that there exists $s_0\geq0$ such that $\mathcal{H}^{h^s}(A)=\infty$
if $s<s_0$ and $\mathcal{H}^{h^s}(A)=0$ if $s>s_0$. This value $s_0$ is called the \emph{Hausdorff dimension} of the set $A$, which we will denote by HD$(A)$.

It was shown by McMullen \cite{mcmullen} that HD$(\mathcal{J}(E_\lambda))=2$ for all $\lambda$, whereas Julia sets of hyperbolic
exponentials have zero Lebesgue measure (the latter result was shown independently by Eremenko and Lyubich, \cite{eremenkolyubich2}).
McMullen further remarked that $\mathcal{J}(E_\lambda)$ always has infinite Hausdorff measure with respect
to the gauge functions $t\mapsto t^2\log^k(1/t)$, where $k\in\mathbb{N}$ is arbitrary.

These results give rise to the question of characterizing the gauge functions $h$ with $\mathcal{H}^h(\mathcal{J}(E_\lambda))=\infty$
(resp. $\mathcal{H}^h(\mathcal{J}(E_\lambda))=0$ if $E_\lambda$ is hyperbolic).

If $0<\lambda<1/e$, $E_\lambda$ has exactly two real fixed points $\alpha_\lambda$ and $\beta_\lambda$, where $\alpha_\lambda$ is \emph{attracting} (i.e. $E_\lambda'(\alpha_\lambda)<1$) and $\beta_\lambda$ is \emph{repelling} (i.e. $E_\lambda'(\beta_\lambda)>1$). Recall that a classical result
of K\oe nigs implies that there exists a holomorphic function $\Phi_\lambda$, defined in a neighborhood $V$ of $\beta_\lambda$, which satisfies $\Phi_\lambda(\beta_\lambda)=0$, $\Phi_\lambda'(\beta_\lambda)=1$ and
\begin{equation}\label{99}
\Phi_\lambda(E_\lambda(z))=\beta_\lambda\Phi_\lambda(z)\text{ for all }z\text{ such that
}z,E_\lambda(z)\in V.
\end{equation}
The proof of K\oe nigs' theorem (and of some other results stated here without proof) can be found in
standard monographs about complex dynamics (e.g. \cite{beardon}, \cite{milnor}, \cite{steinmetz}).

It is easy to see that $\Phi_\lambda(U\cap\mathbb{R})\subset\mathbb{R}$ and that \eqref{99} admits
a real-valued continuation of $\Phi_\lambda$ to $[\beta_\lambda,\infty)$. Further, $\Phi_\lambda(x)$ tends
to $\infty$ as $x\to\infty$, but slower than any iterate $\log^k$ of the logarithm.

Our main result is the following:
\begin{theo}\label{main1}
Let $\lambda_0\in(0,1/e)$. Define $\beta_{\lambda_0}$ and $\Phi_{\lambda_0}$ as above.
Let $K_{\lambda_0}:=\frac{\log 2}{\log\beta_{\lambda_0}}$ and let $h(t)=t^2g(t)$ be a gauge function.
\begin{enumerate}
\item If
$$\liminf_{t\to0}\frac{\log
g(t)}{\log\Phi_{\lambda_0}(1/t)}>K_{\lambda_0},$$ then
$\mathcal{H}^h(\mathcal{J}(E_\lambda))=\infty$ for every
$\lambda\in\Com\setminus\set{0}$. The measure $\mathcal{H}^h$ is
not even $\sigma$-finite on $\mathcal{J}(E_\lambda)$.
\item If
$$\limsup_{t\to0}\frac{\log
g(t)}{\log\Phi_{\lambda_0}(1/t)}<K_{\lambda_0},$$ then
$\mathcal{H}^h(\mathcal{J}(E_\lambda))=0$ if $\lambda\in(0,1/e).$
\end{enumerate}
\end{theo}
We will prove this theorem by showing that if $h(t)=h_{\lambda_0,\gamma}(t):=t^2\Phi_{\lambda_0}(1/t)^\gamma$, then
$\mathcal{H}^h(\mathcal{J}(E_\lambda))=\infty$ when $\gamma>\log 2/\log\beta_{\lambda_0}$
and (provided that $\lambda\in(0,1/e)$) $\mathcal{H}^h(\mathcal{J}(E_\lambda))=0$ when $\gamma<\log 2/\log\beta_{\lambda_0}$. Hence we will use the
functions $h_{\lambda,\gamma}$ very frequently.

We show further that statement (a) of the above theorem remains true if we
replace $\mathcal{J}(E_\lambda)$ by $I(E_\lambda)$,
and that statement (b) is still valid if $\lambda$ is a hyperbolic parameter or
if we replace $\mathcal{J}(E_\lambda)$ by $I(E_\lambda)$ and $\lambda$ is
arbitrary. In (b), the constant $K_{\lambda_0}$ then has to be
replaced by a smaller constant which not only depends on $\lambda_0$, but
also on $\lambda$.

There is a vast literature on dynamics of exponential maps. We mention only a few
results. Devaney and Krych
\cite{devaneykrych} proved that $\mathcal{J}(E_\lambda)$ is a 'Cantor bouquet' for $0<\lambda<1/e$, i.e.
it is homeomorphic to the product of a Cantor set and the line $[0,\infty)$. Schleicher and Zimmer \cite{schleicherzimmer} proved that
$I(E_\lambda)$ consists of curves for all $\lambda$. Karpi\'{n}ska (\cite{karpinska2},\cite{karpinska})
proved the following dimension paradox: For $0<\lambda<1/e$, the Hausdorff dimension of the
endpoints of the curves which form the Cantor bouquet is 2, but the Hausdorff dimension
of the curves without endpoints is 1. Urba\'{n}ski and Zdunik \cite{urbanskizdunik} showed
that the Hausdorff dimension of the set of non-escaping points in the Julia set of a hyperbolic
exponential map is always less than 2. Further, Rempe \cite{rempe} proved an analogue of Böttcher's theorem for exponential maps.
The last two results mentioned will be used later in this work.

For an introduction to the dynamics of exponential maps, we refer the reader to the extensive surveys by Devaney \cite{devaney} 
and Schleicher \cite{schleicher}, as well as Rempe's article
mentioned above. For an introduction to general transcendental dynamics, see e.g. Bergweiler \cite{bergweiler}.

\subsection{Outline of the paper}
In section 2, we provide some notations that we will use throughout this work. After that, we review some
basic results from function theory and apply them
to obtain results about distortion of holomorphic maps.

In section 3, we develop a sufficient condition for a gauge
function $h$ such that $\mathcal{J}(E_\lambda)$ has infinite (not
even $\sigma$-finite) Hausdorff measure with respect to $h$.
The ideas of the proof are due to McMullen, but we have to estimate things more
carefully.

Section 4 gives a sufficient condition for a gauge function $h$
such that $\mathcal{J}(E_\lambda)$ has zero Hausdorff measure with
respect to $h$. In contrast to section 3, this $h$ depends on $\lambda$. Further, $\lambda$ is restricted to values between
0 and $1/e$.

In section 5, we show that the restriction that $h$ depends on $\lambda$ in unnecessary by showing that if
$\beta_{\lambda_1}^{\gamma_1}=\beta_{\lambda_2}^{\gamma_2}$, the
two resulting gauge functions from section 4 have the same growth,
i.e. their quotient can be estimated from above and below by
positive constants. After that, we prove
Theorem \ref{main1} and describe how its statement can be generalized to Julia sets
of hyperbolic exponential maps and escaping sets of arbitrary exponential maps.

\subsection*{Acknowledgements}
I deeply thank Walter Bergweiler for his constant support and advice. I am grateful to Albert Clop, Adam Epstein,
Boguslawa Karpi\'{n}ska, Janina Kotus, Lasse Rempe and Dierk Schleicher
for many interesting and fruitful discussions about this work.
Further, I thank the EU Research Training Network CODY for their financial support.

\section{Notations and preliminaries}

\subsection{Notations}

For $z \in \Com$, let $\Re z$ and $\Im
z$ denote the real and imaginary parts of $z$. If $z_0 \in \Com$
and $r>0$, we write $D(z_0,r)$ for the disc in $\Com$ with center
$z_0$ and radius $r$ with respect to the euclidean metric. By
$\mathbb{D}:=D(0,1)$ we denote the open unit disc in $\Com$. For
$\theta \in \mathbb{R}$, let $Q(z_0,r,\theta)$ be the square with center $z_0$ and side
length $r$ whose sides have angle $\theta$ with the coordinate axes. Thus
$$Q(0,r,0)=\set{z\in\Com\Trenn\max\set{\Re z,\Im z}<\frac{r}{2}}$$
and $$Q(z_0,r,\theta)=z_0+e^{i\theta}Q(0,r,0).$$ If the angle
$\theta$ is not important, we will suppress it and just write
$Q(z_0,r)$ in order to increase readability.
If $A \subset \Com$ is Lebesgue-measurable, we denote its Lebesgue
measure by $\abs{A}$. If $A,B \subset \Com$ are (Lebesgue-)
measurable and $0<\abs{B}<\infty$, we write dens($A,B$) for the
density of $A$ in $B$, i.e.
$$\text{dens}(A,B):=\frac{\abs{A \cap B}}{\abs{B}}.$$ If $x$ is a real number,
we denote by $\left[ x\right]$ the largest integer which is not greater than $x$.
We denote the \emph{postcritical set} of $E_\lambda$ by $P(E_\lambda)$, i.e.
$$P(E_\lambda):=\overline{\bigcup_{n\in\mathbb{N}_0}E_\lambda^n(0)}.$$
For the remainder of this paper, let $\lambda'\in(0,1/e)$ be fixed. Set $E:=E_{\lambda'}$, $\beta:=\beta_{\lambda'}$, $\Phi:=\Phi_{\lambda'}$ and $h_\gamma:=h_{\lambda',\gamma}$ for $\gamma>0$ in order to suppress indices (for
the definition of the right sides, see section 1).

\subsection{Basic results}

In this section, we provide most of the tools that are needed to prove Theorem \ref{main1}.

We start with the so-called 'blow-up property' of Julia sets. Recall that the backward orbit
of $z\in\Com$ consists of all points $w\in\Com$ such that $f^n(w)=z$ for some $n\in\mathbb{N}$.
\begin{lemm}\label{blow}
Let $f$ be an entire function.
If $U\subset\Com$ is an open set intersecting $\mathcal{J}(f)$ and
$K\subset\Com$ is a compact set which does not contain any point with finite backward orbit, then there exists $n_0\in\mathbb{N}$ such that
$f^n(U)\supset K$ for all $n\geq n_0$.
\end{lemm}
The next result follows directly from Cauchy's integral formula:
\begin{lemm}\label{cif}
If a holomorhic function $f$ maps $D(z_0,r)$ into a disk of radius
$s$, then
$$\abs{f'(z_0)}\leq\frac{s}{r}.$$
\end{lemm}
The following lemma is a simple application of the Koebe growth and distortion theorems.
\begin{lemm} \label{kk}
Let $z_0 \in \Com$, $r>0$, $f:D(z_0,r) \to \Com$ be a univalent
function and $z \in D(z_0,r)$. Then
\begin{align} \label{9}&r^2\abs{f'(z_0)}\frac{r-\abs{z-z_0}}{(r+\abs{z-z_0})^3} & \leq
& \hspace{1.1cm} \abs{f'(z)} && \leq &
r^2\abs{f'(z_0)}\frac{r+\abs{z-z_0}}{(r-\abs{z-z_0})^3}
\\ \text{and} \notag \\ \label{10}
&r^2\abs{f'(z_0)}\frac{\abs{z-z_0}}{(r+\abs{z-z_0})^2} & \leq &
\hspace{0.4cm} \abs{f(z)-f(z_0)} && \leq &
r^2\abs{f'(z_0)}\frac{\abs{z-z_0}}{(r-\abs{z-z_0})^2}
\end{align}
\end{lemm}
Koebe's theorems imply in particular that the class $S$ of all univalent functions $f:\mathbb{D}\to\Com$
such that $f(0)=0$ and $f'(0)=1$ is normal. This yields the following result.
\begin{theo} \label{epsdelta}
For every $\eps>0$, there exists $\delta>0$ such that if $f \in S$ and $z \in
\mathbb{D}$ with $\abs{z}<\delta$, then
$$\abs{\frac{f(z)}{z}-1}<\eps.$$
\end{theo}
Let $A\subset\Com$ either be open and bounded or the closure of such a set.
A function $f:A\to\Com$ is said to have \emph{bounded distortion} if $f$ is a bilipschitz mapping, i.e.
$$0<c_f:=\inf_{\substack{z,w \in A\\z\neq w}}\frac{\abs{f(z)-f(w)}}{\abs{z-w}}\leq
\sup_{\substack{z,w \in A\\z\neq
w}}\frac{\abs{f(z)-f(w)}}{\abs{z-w}}=:C_f<\infty.$$
The distortion of $f$ is then defined as $D(f):=C_f/c_f$.

It can be easily shown that if $A$ is open and $f$ has bounded distortion, then $f$ extends to a function on $\overline{A}$
with the same distortion as $f$. Conversely, if $A$ is the closure of an open bounded set, then $f|_{\text{int}(A)}$ has the same
distortion as $f$. Note that maps with bounded distortion are injective.

Some simple properties of holomorphic functions with bounded distortion are summarized in the following lemma.
\begin{lemm} \label{eigensch}
Let $U \subset \Com$ be bounded and open. Let $f:U \to \Com$ be a holomorphic
function with bounded distortion. Then the following statements
hold:
\begin{enumerate}
\item $D(f)=D(f^{-1})$ \\
\item If $V\supset f(U)$ is a domain and $g:V\to\Com$ is
univalent with bounded distortion, then $D(g\circ f)\leq
D(g)D(f)$ \\
\item If $A\subset U$ is Lebesgue-measurable, then $f(A)$ is
Lebesgue-measurable and
$$dens(A,U)\leq D(f)^2dens(f(A),f(U))$$
\item For every $z \in U$, $$c_f\leq\abs{f'(z)}\leq c_fD(f)$$
\end{enumerate}
\end{lemm}
Using the Koebe theorems, we obtain a simple result which will allow us to estimate
the distortion of a map:
\begin{lemm} \label{first}
Let $z_0 \in \Com,r>0,K>3$. Let $f:D(z_0,Kr)\to\Com$ be a
univalent function. Then
$$D\left(f|_{D(z_0,r)}\right)<\Big(\frac{K+1}{K-3}\Big)^6.$$
\end{lemm}
From this result we obtain immediately
\begin{lemm}\label{applic}
Let $Q_1,\ldots,Q_n$ be a finite sequence of squares in $\Com$
with the same side length $r$, i.e. $Q_i=Q(a_i,r)$ for some $a_i\in\Com$. Let $K>3$ and
$f_i:Q_i \to \Com$ be univalent maps such that $Q_{i+1}\subset
f_i(Q_i)$ for every $1\leq i\leq n-1$. Suppose that both $f_n$ and
$f_1^{-1}\circ\ldots\circ f_{n-1}^{-1}|_{Q_n}$ can be extended
univalently to $D(a_n,Kr/\sqrt{2})$. Let
$$F:=(f_n\circ\ldots\circ f_1)^{-1}:f_n(Q_n)\to Q_1.$$ Then
$$D(F)\leq\Big(\frac{K+1}{K-3}\Big)^{12}$$ independent of $n$.
\end{lemm}
In particular, if there exists an entire function $f$ such that all the $f_i$ are restrictions of
$f$ to $Q_i$, then $(f^n)^{-1}$ (restricted to $f(Q_n)$) has uniformly
bounded distortion. This was already shown by McMullen \cite{mcmullen} using different techniques.

Now we are ready to state the result that any square which is small enough will be
mapped under an arbitrary univalent function $f$ to a set which is almost square-shaped. This result is
an application of Theorem \ref{epsdelta}.
\begin{lemm}\label{second}
For every $1/\sqrt{2}>\eps>0$, there exists a constant $K>1$ with
the following property: Let $z_0 \in
\Com$, $r>0$, $K'\geq K$, a univalent function
$\tilde{f}:D\left(z_0,K'r/\sqrt{2}\right) \to \Com$ and a
square $Q=Q(z_0,r,\theta)$ be given. Let
$f:=\tilde{f}|_{\overline{Q}}$ and $d:=D(f)$ be the distortion of
$f$. Then
$$Q\left(f(z_0),\abs{f'(z_0)}r\frac{1}{d}\bra{1-\sqrt{2}\eps},\theta+\arg
f'(z_0)\right)\subset f(Q)$$ and $$f(Q)\subset
Q\left(f(z_0),\abs{f'(z_0)}rd\bra{1+\sqrt{2}\eps},\theta+\arg
f'(z_0)\right).$$
\end{lemm}

\section{The estimate from below}

In this section, we prove that if $\beta^\gamma>2$,
then the Julia set of any exponential map $E_\lambda$ has infinite (even non-$\sigma$-finite) Hausdorff measure
with respect to $h_\gamma$.

\subsection{Preparations}

We keep this section rather short, since most of the
work done here is an obvious generalization of methods developed by McMullen or
consists of simple calculations.

We first define what it means for a family of sets to satisfy the 'nesting conditions'.
\begin{defi}[nesting conditions]\label{nest}
For $k \in \mathbb{N}$, let $\mathcal{A}_k$ be a finite collection
of compact, disjoint and connected subsets of $\Com$ with positive
Lebesgue-measure. Let $A_k$ be the union of the elements of $\mathcal{A}_k$. We
say that $\set{\mathcal{A}_k}$ satisfies the \emph{nesting conditions} if it has the following three properties:
\begin{enumerate}
\item For every $k \in \mathbb{N}$ and $F \in \mathcal{A}_{k+1}$, there exists some $F' \in
\mathcal{A}_k$ such that $F \subset F'$.
\item There exists a decreasing sequence $(d_k)$
converging to 0 such that $$\max_{F\in\mathcal{A}_k}\text{diam}(F)\leq d_k \text{ for all }k\in\mathbb{N}$$
\item There exists a sequence $(\Delta_k)$ of positive reals such
that $$\text{dens}(A_{k+1},F)\geq\Delta_k\text{ for all
}k\in\mathbb{N},F\in\mathcal{A}_k.$$
\end{enumerate}
The intersection $A:=\bigcap_{k=1}^\infty A_k$ is a non-empty and compact set.
\end{defi}
We will use the following Frostman-type lemma (see e.g. \cite{urbanski}, Theorem 7.6.1).
\begin{lemm} \label{frostman} Let $h$ be a continuous gauge function and
$\mu$ be a Borel probability measure on $\Com$.
Let $X \subset \Com$ be $\mu$-measurable and $0<c<\infty$. If
\begin{equation} \label{5}
\limsup_{r \to 0} \frac{\mu(D(z,r))}{h(r)}\leq c \text{ for all }z \in
X, \end{equation} then $\mathcal{H}^h(X)\geq\mu(X)/c$.
\end{lemm}
The key lemma to the proof of this section's main result is
\begin{lemm} \label{infinity}
Let $\set{\mathcal{A}_k}$ be a collection of families of sets
which satisfies the nesting conditions $($with properly chosen
sequences $(d_k)$ and $(\Delta_k))$. Let $A$ be defined as above.
Let $\eps>0$ and $g:(0,\eps)\to\mathbb{R}_{\geq 0}$ be a
decreasing continuous function such that $t^2g(t)$ is increasing.
Further, suppose that
$\lim_{t\to 0}t^2g(t)=0$ and
\begin{equation}\label{69}
\lim_{k\to\infty} g(d_k)\prod_{j=1}^k\Delta_j=\infty.
\end{equation}
Define
$$h:[0,\eps)\to\mathbb{R},t\mapsto
\begin{cases}
t^2g(t)&,t>0\\
0&,t=0
\end{cases}$$
Then $h$ is a continuous gauge function and we have
$\mathcal{H}^h(A)=\infty$.
\end{lemm}
\begin{proo}
The proof follows ideas of McMullen (\cite{mcmullen}, Proposition 2.2), we only give
a sketch here. First note that $h$ is clearly a continuous gauge function.

For $F\in \mathcal{A}_1$, we define
$\tau_1(F):=\abs{F}/\abs{A_1}.$ For
$F\in\mathcal{A}_{k+1}$, there exists a unique $G\in\mathcal{A}_k$
such that $F\subset G$. Denoting by $F_i$ the elements of
$\mathcal{A}_{k+1}$ that are contained in $G$, we define
$$\tau_{k+1}(F):=\frac{\abs{F}}{\sum\limits_i\abs{F_i}}\tau_k(G).$$
Now we set $\mathcal{X}:=\bigcup_{k \in
\mathbb{N}}\mathcal{A}_k\cup\bigcup_{k\in\mathbb{N}}\set{M\subset\Com\trenn
M\cap A_k=\emptyset}$
and define $\tau$ by
$$\tau:\mathcal{X}\to \mathbb{R}_{\geq 0}, B \mapsto
\begin{cases}
\tau_k(B)&,B\in\mathcal{A}_k\\
0&,B\notin\bigcup_{k\in\mathbb{N}}\mathcal{A}_k
\end{cases}$$
By setting $$\eta(B):=\lim_{\delta\to 0}\inf\set{\sum_{i=1}^\infty \tau(B_i)\Trenn
B\subset\bigcup_{i\in\mathbb{N}}B_i,B_i\in
\mathcal{X},\text{diam}(B_i)<\delta},$$
we obtain a metric outer measure on $\Com$.
The restriction $\mu$ of $\eta$ to all $\eta$-measurable sets
is a Borel probability measure supported on $A$ with $\mu(A)=1$. It is easy to
see that $\mu$ coincides with $\tau$ on every $\mathcal{A}_k$.

Let $z\in A$, $r>0$ and $D:=D(z,r)$. Choose $k\in\mathbb{N}$ such that
$d_k>r\geq d_{k+1}$. Note that $k$ tends to infinity if and only if $r$ tends to 0.
Let $\tilde{D}$ be the union of all sets in $\mathcal{A}_{k+1}$
which intersect $D$. Then $\text{diam}(\tilde{D})\leq 4r$. Since $\mu$ is supported on $A$,
we see that $\mu(D)\leq\mu(\tilde{D})+\mu(D\setminus\tilde{D})=\mu(\tilde{D})$.
Further, it is not difficult to show that $\mu(\tilde{D})\leq\frac{|\tilde{D}|}{\abs{A_1}}\frac{1}{\Delta_1\cdots\Delta_k}$.
Combining the previous three estimates yields $\mu(D)\leq\frac{\abs{D(z,4r)}}{\abs{A_1}}\frac{1}{\Delta_1\cdots\Delta_k}$.
Since
$$\frac{\abs{D(z,4r)}}{\abs{A_1}}\frac{1}{\Delta_1\cdots\Delta_k}=Cr^2g(r)\frac{1}{g(r)
\Delta_1\cdots\Delta_k}\leq C\frac{1}{g(d_k)\Delta_1\cdots\Delta_k}h(r)$$
and $C/(g(d_k)\Delta_1\cdots\Delta_k)\to 0$ as $k\to\infty$ by the hypothesis,
an application of Lemma \ref{frostman} finishes the proof.
\end{proo}

Our goal is to construct a family of sets that, given $\eps>0$, satisfies the nesting conditions with $\Delta_k>\frac{1}{2}-\eps$, $d_k<1/E^k(2\beta)$ and $A\subset\mathcal{J}(E_\lambda)$. Then we can use the above lemma
and the functional equation for $\Phi$ to obtain the desired result.

In the remaining part of this section, let $\lambda\in\Com$ be fixed. For $\delta>0$, we define
$$I_\delta:=\Big[-\frac{\pi}{2}+\delta-\arg\lambda,\frac{\pi}{2}-\delta-\arg\lambda\Big]$$
and $$J_\delta:=\set{z\in\Com\trenn\Im z\in (I_\delta+2\pi ik)\text{
for some }k\in\mathbb{Z}}.$$ If $z=x+iy\in J_\delta$, then
$\arg E_\lambda(z)=\arg\lambda+y$, and hence
\begin{equation}\label{42}
\Re E_\lambda(z)=\abs{\lambda}e^{\Re z}\cos\arg
E_\lambda(z)\geq\abs{\lambda}\cos\Big(\frac{\pi}{2}-\delta\Big)e^{\Re
z}=:Ce^{\Re z}.
\end{equation}
Consequently the real part of
$E_\lambda^n(z)$ grows exponentially fast as long as the iterates $E_\lambda^n(z)$ stay in
$J_\delta$ and $\Re z$ is large enough, so $z\in
I(E_\lambda)\subset\mathcal{J}(E_\lambda)$ if
\begin{itemize}
\item $E_\lambda^n(z)\in J_\delta$ for all
$n\in\mathbb{N}$ and
\item $z$ belongs to a suitable right half plane.
\end{itemize}
To estimate the densities $\Delta_k$ of our sets later on, we need the following
lemma, which we do not prove it in full detail since its statement is quite obvious.
We introduce a notation first and state the lemma afterwards.
\begin{defi}[$r$-box, $r$-packing]
For $r>0$, a subset $Q\subset\Com$ is an \emph{$r$-box} if
$Q=\overline{Q(z_0,r,0)}$ for some $z_0\in\Com$, i.e. a square of
side length $r$ with sides parallel to the coordinate axes. For a
set $B\subset\Com$ and $r>0$, an \emph{$r$-packing} of $B$ is a union
of disjoint $r$-boxes which are contained in $B$.
\end{defi}
\begin{lemm}\label{pack}
For every $\delta,\delta'>0$, there exists $r_0>0$ with the
following property: If $r<r_0$, we can find $c(r)>0$ such
that for every $z_0\in\Com,c>c(r)$ and $\theta\in\mathbb{R}$, the
set $Q:=Q(z_0,cr,\theta)\cap J_\delta$ has an $r$-packing
\textnormal{pack}$_Q$ which satisfies
$$\textnormal{dens}(\textnormal{pack}_Q,Q)>\frac{1}{2}-\frac{\delta}{\pi}-\delta'.$$
\end{lemm}
\begin{proo}
It is easy to see that for every $r>0$, there exists $c>0$ which satisfies
$$\text{dens}(J_\delta, Q(z_0,rc,\theta))>\frac{1}{2}-\frac{\delta}{\pi}
-\frac{\delta'}{2}\text{ for every
}z_0\in\Com,\theta\in\mathbb{R}.$$ Let $Q=Q(z_0,rc,\theta)$ be an
arbitrary square in $\Com$. In order to find our $r$-packing of $Q$, we fix some
small $\eps>0$ and cover $\Com$ with an $(r+\eps)$-grid of
$(r+\eps)$-boxes. We define an $r$-packing pack$_Q$ of $Q\cap
J_\delta$ as follows: An $r$-box $Q(z,r,0)$ is contained in
pack$_Q$ if and only if $Q(z,r+\eps,0)$ is completely contained in $J_\delta\cap
Q$ and belongs to the
$(r+\eps)$-grid. \\ We want to estimate the number of $(r+\eps)$-boxes which are
contained in pack$_Q$. First, note that at least
$\left[\frac{(rc)^2}{(r+\eps)^2}\right]$ squares from the grid
intersect $Q$, because otherwise the sum of the areas of the grid
squares would be less than $(rc)^2=\abs{Q}$ which is a
contradiction. Since
$$\abs{J_\delta\cap
Q}>\Big(\frac{1}{2}-\frac{\delta}{\pi}-\frac{\delta'}{2}\Big)\abs{Q},$$
it follows that the number of grid squares that intersect
$J_\delta\cap Q$ is at least
$$\left[\Big(\frac{1}{2}-\frac{\delta}{\pi}-\frac{\delta'}{2}\Big)\frac{(rc)^2}{(r+\eps)^2}\right].$$
The number of grid squares that intersect $Q$ and the boundary of
a particular strip $T$ in $J_\delta$ is bounded by
$\mathcal{O}\big(\frac{rc}{r+\eps}\big)=\mathcal{O}(c)$, and there
are $\mathcal{O}(rc)$ strips in $J_\delta$ intersecting $Q$. Hence
it follows that the number of grid squares that intersect
$\partial J_\delta\cap Q$ is bounded by $\mathcal{O}(rc^2)$. Consequently,
\begin{align*}
\text{dens}(\text{pack}_Q,Q)
&\geq\frac{(r+\eps)^2}{r^2c^2}\bigg(\Big(\frac{1}{2}-\frac{\delta}{\pi}-\frac{\delta'}{2}\Big)\frac{r^2c^2}{(r+\eps)^2}-1-\mathcal{O}(rc^2)\bigg)\\
&=\frac{1}{2}-\frac{\delta}{\pi}-\frac{\delta'}{2}+\mathcal{O}\Big(\frac{1}{c^2}\Big)+\mathcal{O}(r)\\
&>\frac{1}{2}-\frac{\delta}{\pi}-\delta'
\end{align*}
if $r$ is small enough and $c$ is large enough.
\end{proo}

The next result, which is mainly an application of the preceding
lemma, will help us to estimate the densities of our nested sets.
\begin{lemm}\label{dens}
For every $\tilde{\eps}>0$, there exist $\delta>0$ and $r_0>0$
with the following property: If $r\leq r_0$, we can find $x>0$
such that if $$Q:=Q(z_0,r,\theta)\subset(J_\delta\cap
\set{z\trenn\Re z>x}),$$ then there exists an r-packing
$\textnormal{pack}_{E_\lambda(Q)}$ of $E_\lambda(Q)\cap J_\delta$
with
$$\textnormal{dens}(\textnormal{pack}_{E_\lambda(Q)},E_\lambda(Q))>\frac{1}{2}-\tilde{\eps}.$$
\end{lemm}
\begin{proo}
Let $\eps>0$. Because $E_\lambda$ is univalent on $D(z,\pi)$ for
all $z\in\Com$, we conclude from Lemma \ref{second} that there exists $r_1>0$ such that for every $r'<r_1$
and every $Q=Q(z,r',\theta)$, we can find a square
$$Q\Big(E_\lambda(z),\abs{E_\lambda'(z)}r'\frac{1}{D(E_\lambda|_Q)}(1-\sqrt{2}\eps)\Big)$$
which is contained in $E_\lambda(Q)$ and another square
$$Q\Big(E_\lambda(z),\abs{E_\lambda'(z)}r'D(E_\lambda|_Q)(1+\sqrt{2}\eps)\Big)$$
that contains $E_\lambda(Q)$. By Lemma \ref{first}, we can choose $r_2$ so small
that $D(E_\lambda|_Q)<1+\eps$ for every $Q=Q(z,r',\theta)$,
$r'<r_2$. Hence
$$Q\Big(E_\lambda(z),\abs{E_\lambda(z)}r'\frac{1}{1+\eps}(1-\sqrt{2}\eps)\Big)\subset
E_\lambda(Q)\subset
Q\Big(E_\lambda(z),\abs{E_\lambda(z)}r'(1+\eps)(1+\sqrt{2}\eps)\Big)$$
for all $r'<\min\set{r_1,r_2}$ and all $Q=Q(z,r',\theta)$. Now choose $\delta,\delta'>0$ such that
$$\frac{1}{2}-\frac{\delta}{\pi}-\delta'>\frac{1}{2}-\eps.$$ By Lemma \ref{pack}, there
exists $r_3>0$ with the following property: \\ If $r'<r_3$, we can find $c(r')>0$ such that for every $Q=Q(z,cr',\theta)$ (where $c\geq c(r')$),
the set $Q\cap\mathcal{J}_\delta$ has an $r'$-packing pack$_Q$ with
\begin{equation}\label{49}
\text{dens}(\text{pack}_Q,Q)>\frac{1}{2}-\eps.
\end{equation}
So choose $r<\min\set{r_1,r_2,r_3}$ and $c(r)>0$ with the above
property. Then we can find $x>0$ such that
$$\abs{E_\lambda(z)}>c(r)\frac{1+\eps}{1-\sqrt{2}\eps}$$ for every
$z\in\Com$ with $\Re z>x$. If we fix some square
$Q=Q(z_0,r,\theta)\subset\set{z\trenn\Re z>x}$, it follows that
\begin{align*}
Q_1:=Q\Big(E_\lambda(z_0),\abs{E_\lambda(z)}r\frac{1-\sqrt{2}\eps}{1+\eps}\Big)&\subset
E_\lambda(Q)\\
&\subset
Q\Big(E_\lambda(z_0),\abs{E_\lambda(z)}r(1+\eps)(1+\sqrt{2}\eps)\Big)=:Q_2.
\end{align*}
Let an $r$-packing pack$_{Q_1}$ of $Q_1\cap J_\delta$ which
satisfies \eqref{49} (with $Q_1$ instead of $Q$) be given. Then
pack$_{Q_1}$ is also an $r$-packing of $E_\lambda(Q)\cap
J_\delta$. Thus
\begin{align*}
\text{dens}(\text{pack}_{Q_1},E_\lambda(Q))&=
\frac{\abs{\text{pack}_{Q_1}\cap
E_\lambda(Q)}}{\abs{E_\lambda(Q)}}\\&\geq\frac{\abs{\text{pack}_{Q_1}\cap
Q_1}}{\abs{Q_2}}\\&=\frac{\abs{\text{pack}_{Q_1}\cap
Q_1}}{\abs{Q_1}}\frac{\abs{Q_1}}{\abs{Q_2}}\\&>\Big(\frac{1}{2}-\eps\Big)\frac{(1-\sqrt{2}\eps)^2}{(1+\sqrt{2}\eps)^2(1+\eps)^4}.
\end{align*}
Because $\eps>0$ was arbitrary, we can choose $\eps$ so small that
the last term is smaller than $\frac{1}{2}-\tilde{\eps}$.
\end{proo}

If we want to estimate the distortion of inverse branches of the
exponential map with the help of Lemma \ref{applic}, we have to be
sure that these inverse branches are defined on a large
region. This is ensured by the following obvious lemma.
\begin{lemm}\label{psaway}
Let $\delta>0$ and $J_\delta$ be defined as above. Let $K>0$. Then
there exists $x_0>0$ with the following property:\\ If $z\in\Com,n\in\mathbb{N}_0$
with $\Re z>x_0$ and $E_\lambda^k(z)\in J_\delta$ for all $0\leq
k\leq n$, then
$$\Re E_\lambda^n(z)-\Re E_\lambda^k(0)\geq K\text{ for all }0\leq k\leq n.$$
\end{lemm}
Finally, the following result will be used to prove the non-$\sigma$-finiteness. Its proof can for example be found in \cite{rogers}.
\begin{lemm}\label{nsigma}
Let $h_1$ and $h_2$ be gauge functions with $\frac{h_2(t)}{h_1(t)} \to 0 \text{ as } t \to 0$ and $A
\subset \Com$. If $A$ has $\sigma$-finite
$\mathcal{H}^{h_1}$-measure, then $A$ has zero
$\mathcal{H}^{h_2}$-measure.
\end{lemm}

\subsection{A gauge function which leads to infinite measure}

After these preparations, we can prove the main result of this section.

\begin{theo}\label{infmeas}
Let $\lambda\in\Com\setminus\set{0}$ and $\gamma>0$ such that $\beta^\gamma>2$.
Then $\mathcal{H}^{h_\gamma}(\mathcal{J}(E_\lambda))=\infty.$
\end{theo}
\begin{proo}
First of all, let $\eps>0$ such that
\begin{equation}\label{60}
\beta^\gamma\frac{1/2-\eps}{(1+\eps)^2}>1.
\end{equation}
Further, let $K>0$ be so large that
\begin{equation}\label{64}
\Big(\frac{K+1}{K-3}\Big)^{12}<1+\eps.
\end{equation}
By Lemma \ref{dens}, there exist $\delta>0$, $x_0>0$ and $r>0$
with the following three properties:
\begin{enumerate}
\item
$Kr/\sqrt{2}<\pi$
\item If $Q$ is an $r$-box with $Q\subset J_\delta\cap\set{\Re
z>x_0}$, then there exists an $r$-packing pack$_{E_\lambda(Q)}$ of
$E_\lambda(Q)\cap J_\delta$ such that
\begin{equation*}
\text{dens}(\text{pack}_{E_\lambda(Q)},E_\lambda(Q))>\frac{1}{2}-\eps.
\end{equation*}
\item
$r\sqrt{2}\lambda'/\abs{\lambda}<1$
\end{enumerate}
For every $Q$ as in (b), let $\mathcal{B}(E_\lambda(Q))$ denote
the collection of the $r$-boxes that form pack$_{E_\lambda(Q)}$.
Choose $x_1\geq x_0$ such that
\begin{equation}\label{67}
\Re E_\lambda^n(z)>E^n(2\beta)
\end{equation}
for all $z\in J_\delta\cap\set{\Re z>x_1}$ and all $n\in\mathbb{N}$
with $E_\lambda^k(z)\in J_\delta$ whenever $0\leq k\leq n$. By Lemma \ref{psaway}, we can choose $x_2\geq x_1$ such that for
all $z\in J_\delta\cap\set{\Re z>x_2}$ and $n\in\mathbb{N}_0$ with
$E_\lambda^k(z)\in J_\delta$ for $0\leq k\leq n$, we have
\begin{equation}\label{63}
\Re E_\lambda^n(z)-\Re E_\lambda^k(0)\geq\frac{Kr}{\sqrt{2}}\text{
whenever }0\leq k\leq n.
\end{equation}
We define $T:=J_\delta\cap\set{z\trenn\Re z>x_2}$ and choose an
$r$-box $Q\subset T$. Now we define a family of nested sets as
follows: Set $\mathcal{A}_0:=\set{Q}$ and
$$\mathcal{A}_k:=\set{G\trenn G\subset
F\text{ for some }F\in\mathcal{A}_{k-1}\text{ and
}E_\lambda^k(G)\in\mathcal{B}(E_\lambda^k(F))}$$ for
$k\in\mathbb{N}$, i.e. $\mathcal{A}_1$ consists of all preimages
(under $E_\lambda$) of the $r$-boxes that are contained in
$\mathcal{B}(E_\lambda(Q))$, $\mathcal{A}_2$ consists of those
subsets of elements $F\in\mathcal{A}_1$ which are mapped onto an
$r$-box in $\mathcal{B}(E_\lambda^2(F))$ under $E_\lambda^2$ and
so on. Then $\set{\mathcal{A}_k}$ is a family of nested sets. As
in Definition \ref{nest}, we define $A_k:=\bigcup_{F\in\mathcal{A}_k}F$ and
$A:=\bigcap_{k\in\mathbb{N}}A_k$. Then $z\in A$ implies that
$$E_\lambda^n(z)\in J_\delta\text{ for all }n\in\mathbb{N}_0.$$
It follows from \eqref{67} that $\Re E_\lambda^n(z)\to\infty$ as
$n\to\infty$, and consequently
$z~\in~\mathcal{J}(E_\lambda)$.

We now want to estimate the densities and diameters of the sets in
$\mathcal{A}_k$. If $k\in\mathbb{N}$ and $F\in\mathcal{A}_{k-1}$,
then $E_\lambda^{k-1}(F)$ is an $r$-box $Q=\overline{Q(z_Q,r,0)}$.
First, we estimate the distortion of
$(E_\lambda^k)^{-1}|_{E_\lambda(Q)}:E_\lambda(Q)\to F$. Because of
condition (a), $E_\lambda$ is univalent on $D(z_Q,Kr/\sqrt{2})$.
Further, the branch $\varphi$ of $(E_\lambda^{k-1})^{-1}|_Q$ with
$\varphi(z_Q)\in F$ can be continued univalently to
$D(z_Q,Kr/\sqrt{2})$ because of \eqref{63}.
It follows from Lemma \ref{applic} and \eqref{64} that
\begin{equation}\label{65}
D\big((E_\lambda^k)^{-1}|_{E_\lambda^k(F)}\big)<1+\eps.
\end{equation}
Thus we obtain by condition (b) and Lemma \ref{eigensch} (c) that
\begin{align*}
\frac{1}{2}-\eps&<\text{dens}\Big(\hspace{-0.5cm}\bigcup_{Q\in\mathcal{B}(E_\lambda^k(F))}\hspace{-0.5cm}Q,E_\lambda^k(F)\Big)\\
&\leq
D((E_\lambda^k)^{-1}|_{E_\lambda^k(F)})^2\text{dens}\Big((E_\lambda^k)^{-1}\big(\hspace{-0.5cm}\bigcup_{Q\in\mathcal{B}(E_\lambda^k(F))}\hspace{-0.5cm}Q\big),F\Big)\\
&\leq(1+\eps)^2\text{dens}\Big(\hspace{-0.5cm}\bigcup_{Q\in\mathcal{B}(E_\lambda^k(F))}\hspace{-0.5cm}(E_\lambda^k)^{-1}(Q),F\Big).
\end{align*}
Since
$$A_k\cap
F=\hspace{-0.5cm}\bigcup_{Q\in\mathcal{B}(E_\lambda^k(F))}\hspace{-0.5cm}(E_\lambda^k)^{-1}(Q)\cap
F,$$ it follows that
\begin{equation}\label{68}
\text{dens}(A_k,F)=\text{dens}\Big(\hspace{-0.5cm}\bigcup_{Q\in\mathcal{B}(E_\lambda^k(F))}\hspace{-0.5cm}(E_\lambda^k)^{-1}(Q),F\Big)\geq\frac{1/2-\eps}{(1+\eps)^2}=:\Delta_k.
\end{equation}
Now we turn to the diameters of the sets in $\mathcal{A}_k$. Let $F\in\mathcal{A}_k$ and $u,v\in F$. Because $E_\lambda^k(F)$ is an
$r$-box and hence convex, it follows from the mean value inequality, \eqref{67} and condition (c) that
\begin{equation}\label{70}
\abs{u-v}<\frac{r\sqrt{2}}{\abs{\lambda}}\frac{1}{e^{E^{k-1}(2\beta)}}=\frac{r\sqrt{2}\lambda'}{\abs{\lambda}}\frac{1}{E^k(2\beta)}
<\frac{1}{E^k(2\beta)}=:d_k.
\end{equation}
So $\set{\mathcal{A}_k}$ satisfies the nesting conditions. To
finish the proof, we have to check that
$t\mapsto(\Phi(1/t))^\gamma$ satisfies the condition
\eqref{69} in Lemma \ref{infinity}. But this is now
straightforward and follows mainly from the functional equation
\eqref{99}, combined with \eqref{68} and \eqref{70}. In fact,
\begin{align*}
\Phi\Big(\frac{1}{d_k}\Big)^\gamma\prod_{j=1}^k\Delta_j
&=\Phi\Big(\frac{1}{d_k}\Big)^\gamma\Big(\frac{1/2-\eps}{(1+\eps)^2}\Big)^k\\
&=\Phi(E^k(2\beta))^\gamma\Big(\frac{1/2-\eps}{(1+\eps)^2}\Big)^k\\
&=\Big(\beta^k\Phi(2\beta)\Big)^\gamma\Big(\frac{1/2-\eps}{(1+\eps)^2}\Big)^k\\
&=\Big(\beta^\gamma\Big(\frac{1/2-\eps}{(1+\eps)^2}\Big)\Big)^k\Phi(2\beta)^\gamma
\end{align*}
which tends to infinity as $k\to\infty$ by \eqref{60}. Thus Lemma
\ref{infinity} shows that $\mathcal{J}(E_\lambda)$ has infinite
$\mathcal{H}^{h_\gamma}$-measure.
\end{proo}

It follows immediately from Lemma \ref{nsigma} that $\mathcal{H}^{h_\gamma}$ is not even
$\sigma$-finite on $\mathcal{J}(E_\lambda)$.
Note that the set $A$ constructed in the proof above belongs to $I(E_\lambda)$, so the statement of
the theorem remains true if we replace $\mathcal{J}(E_\lambda)$ by $I(E_\lambda)$.

\section{The estimate from above}

Recall that $E=E_{\lambda'}$, $\beta=\beta_{\lambda'}$, $\Phi=\Phi_{\lambda'}$ and $h_\gamma=h_{\lambda',\gamma}$, where
$\lambda'\in(0,1/e)$ is fixed.
In this section, we prove conversely that if $\gamma>0$ is chosen such that $\beta^\gamma<2$, then
$\mathcal{J}(E)$ has zero measure with respect to
$h_\gamma$. Note that, in contrast to section 3, the parameter of the exponential map whose Julia set
we measure and the parameter of the gauge function we use for measuring have to be the same. We will
get rid of this restriction in section 5.

\subsection{Preparations}

In \cite{urbanskizdunik}, Urba\'{n}ski and Zdunik proved that the set of non-escaping points in the Julia set of a hyperbolic exponential
map is small:
\begin{theo}\label{uz}
Let $\lambda\in\Com\setminus\set{0}$ such that $E_\lambda$ is
hyperbolic. Then
$$\textnormal{HD}\big(\mathcal{J}(E_\lambda)\setminus
I(E_\lambda)\big)<2.$$
\end{theo}
In particular, this theorem implies that
$\mathcal{J}(E)\setminus I(E)$ has zero
$\mathcal{H}^{h_\gamma}$-measure for every $\gamma>0$,
which means that we can restrict ourselves to $I(E)$ if we
want to prove that $\mathcal{J}(E)$ has zero
$\mathcal{H}^{h_\gamma}$-measure.\\
We show that it is in fact sufficient to consider a suitable
subset of $I(E)$. For this, we need the following
definition:
\begin{defi}\label{ir}
For $\lambda\in\Com\setminus\set{0}$, we define
$$I_R(E_\lambda):=\set{z\in I(E_\lambda)\trenn\Re
E_\lambda^n(z)\geq R\text{ for all }n\in\mathbb{N}_0}.$$
\end{defi}
Note that
\begin{equation}\label{123456}
I(E_\lambda)=\bigcup_{n=0}^\infty
E_\lambda^{-n}(I_R(E_\lambda))
\end{equation}
since $I(E_\lambda)$ is completely invariant.
\begin{lemm}\label{irreicht}
Let $\lambda\in\Com\setminus\set{0}$ and
$\gamma,R>0$. Then
$\mathcal{H}^{h_\gamma}(I_R(E_\lambda))=0$ implies
that $\mathcal{H}^{h_\gamma}(I(E_\lambda))=0$.
\end{lemm}
\begin{proo}
First note that for every $K>0$, there exists $K'>0$ such that $h_\gamma(Kt)\leq K'h_\gamma(t)$.
Hence zero $\mathcal{H}^{h_\gamma}$-measure is preserved by bilipschitz mappings. The statement now
follows easily from \eqref{123456} since $E_\lambda$ is bilipschitz on small discs.
\end{proo}

\subsection{A gauge function which leads to zero measure}

The main theorem of this section is the following:
\begin{theo}\label{zeromeas}
Let $\gamma>0$ such
that $\beta^\gamma<2$. Then $\mathcal{H}^{h_\gamma}(\mathcal{J}(E))=0.$
\end{theo}
The proof of this theorem is rather technical,
but the idea could not be simpler: For large $R$, we cover the set $I_R(E)$ by small
squares of a definite side length. Chosen one of these squares $Q$, we show that every square $Q_0$ inside
which is small enough has zero $\mathcal{H}^{h_\gamma}$-measure. We do this by dividing $Q_0$ into smaller squares
and showing that a certain percentage of these squares does not intersect the Julia set of $E$. By
repeating this method, we construct a sequence of coverings of $I_R(E)$ by sets whose diameters tend to 0,
and an application of the functional equation yields the desired result.
Observe that the idea of the proof originates from the well known (and much easier) result that
a porous set has Hausdorff dimension less than 2.

Before we start proving the theorem, we have to introduce a
notation. In the last section, we were mostly interested in strips
that may contain points in the Julia set (the strips that formed
$J_\delta$). In this section, we are interested in strips which
are completely contained in the Fatou set: For $\delta>0$, we
define $$F_\delta:=\set{z\in\Com\trenn\Im
z\in\Big[\frac{(4k+1)\pi}{2}+\delta,\frac{(4k+3)\pi}{2}-\delta\Big]\hspace{-0.2cm}\mod
2\pi}.$$ Clearly, $F_\delta\subset\mathcal{F}(E)$ if $\delta$ is small. Further it is immediate
that Lemma \ref{dens} remains valid if we replace $J_\delta$ by
$F_\delta$. \\

\begin{proo2}
By Theorem \ref{uz} and Lemma \ref{irreicht}, it suffices to show
that $\mathcal{H}^{h_\gamma}(I_R(E))=0$ for some $R>0$.
Let $\eps'>0$ such that $(1/2+\eps')\beta^\gamma<1$. Let
$\eps>0$ and $M\in\mathbb{N}$ such that
$$(1-\eps)\frac{(1-2\eps)^4}{(1+2\eps)^{18}}\Big(\frac{1-\sqrt{2}\eps}{1+\sqrt{2}\eps}\Big)^4\Big(\frac{M-1}{M}\Big)^2\Big(\frac{1}{2}-\eps\Big)\geq\frac{1}{2}-\eps',$$
$$(1-\eps)\Big(\frac{1-\sqrt{2}\eps}{1+\sqrt{2}\eps}\Big)^4\frac{(1-2\eps)^4}{(1+2\eps)^{16}}\Big(\frac{1}{2}-2\eps\Big)\geq\frac{1}{2}-\eps'$$
and
\begin{equation}\label{12345}
\frac{1}{1+a}\geq\frac{M-1}{M},
\end{equation}
where we define $a$ by
\begin{equation}\label{1234}
1+a:=\frac{(1+\sqrt{2}\eps)(1+2\eps)^4}{(1-\sqrt{2}\eps)(1-2\eps)^2}.
\end{equation}
Now let $K>3$ with
\begin{equation}\label{103}
\Big(\frac{K+1}{K-3}\Big)^6<1+2\eps
\end{equation}
and such that $K$ satisfies the property of Lemma \ref{second}: If
$r'>0,z\in\Com,\theta\in\mathbb{R}$ and
$f:D(z,Kr'/\sqrt{2})\to\Com$ is univalent, then
$f(Q(z,r',\theta))$ is contained in some square
$Q(f(z),\abs{f'(z)}r'(1+\sqrt{2}\eps)d)$ and contains some square
$Q(f(z),\abs{f'(z)}r'(1-\sqrt{2}\eps)/d)$, where
$d=D(f|_{Q(z,r',\theta)})$. \\
Let $r_0>0$ with
\begin{equation}\label{79}
KM\Big(\frac{1+2\eps}{1-2\eps}\Big)^2\frac{r_0}{\sqrt{2}}<\pi
\end{equation}
and such that $r_0$ satisfies the property of Lemma
\ref{pack}:\\ We can find a constant
$A>(\pi+1)\sqrt{\pi}/r_0\sqrt{\eps}$ such that if
$Q=Q(z,Ar_0,\theta)$, then there exists an $r_0$-packing pack$_Q$
of $Q\cap F_\delta$ with
\begin{equation}
\text{dens}(\text{pack}_Q,Q)>\frac{1}{2}-\eps.
\end{equation}
Because of Lemma \ref{dens}, we can also assume that $r_0$ is so
small that we can find $R_1>0$ and $\delta>0$ with the following
property: If $Q:=Q(z_0,r_0,\theta)\subset\set{\Re z>R_1},$ then
$E(Q)\cap F_\delta$ has an $r_0$-packing
pack$_{E(Q)}$ which satisfies
\begin{equation}
\text{dens}\big(\text{pack}_{E(Q)},E(Q)\big)>\frac{1}{2}-\eps.
\end{equation}
Further it is clear that there exists a universal constant $C'$
with the following property: If we partition $\Com$ into a grid of
squares with side lengths $r_-$ and take any square $Q(z,r_+,\theta)$
with $r_+>r_-$, then at most $C'r_+/r_-$ squares of the $r_-$-grid
intersect the boundary of $Q(z,r_+,\theta)$. Now let
$$B:=90(A+1)r_0$$ and choose $C>0$ such that
$$C>Mr_0\sqrt{2}+B$$ and
$N:=(C-Mr_0\sqrt{2})r_0\frac{1-\sqrt{2}\eps}{1+2\eps}$ satisfies
$$1-\frac{C'}{N}>1-\eps.$$ In particular, we have $C>4KMr_0$. Let
$$L:=C-Mr_0\sqrt{2}>B$$ and $$\kappa:=\frac{1}{LB}.$$
Finally, let $R\geq R_1+(A+1)r_0\sqrt{2}$ such that
\begin{equation}\label{78}
\lambda'e^{R-(A+1)r_0\sqrt{2}}-\alpha\geq\text{dist}(\set{\Re
z\geq R},P(E))\geq C+\sqrt{2}Mr_0,
\end{equation}
where $\alpha$ denotes the real attracting fixed point of $E$.
After these preparations, choose some $z_0'\in\Com$ such that
$$Q_0:=Q(z_0',r_0,0)\cap I_R(E)\neq\emptyset,$$ and
some $x=x(Q_0)\in\mathbb{R}$ with the following five properties:
\begin{enumerate}
\item $x\geq L\, (E^2)'(\max\set{\Re z\trenn z\in Q_0})$
\item $-\frac{\log\kappa}{\log x}\leq\sqrt{\kappa}-\kappa$
\item $E^n(2x)\geq(E^n)'(x)\text{ for all
}n\in\mathbb{N}$
\item $1/E(x)<r_0$
\item
$(E^{n+1})'(x)+\sqrt{2}(E^n)'(x)<\sqrt{2}(E^{n+1})'(x)\text{
for all }n\in\mathbb{N}$
\end{enumerate}
Note that $x$ depends on $Q_0$ only by property (a). From
the first condition, it follows immediately that
\begin{equation}\label{87}
(E^m)'(x)\geq L\abs{(E^{m+2})'(z)}\text{ for all
}m\in\mathbb{N}_0\text{ and }z\in Q_0.
\end{equation}
We define inductively a sequence of side
lengths $r_n$ as follows. Let $$r_1:=\frac{1}{E(x)}$$ and for
$n\in\mathbb{N}$, let $r_{n+1}$ be such that
$$\frac{r_n}{r_{n+1}}\in\mathbb{N}\text{ and
}r_n(E^{n+1})'(x)\in\Big(\frac{r_n}{r_{n+1}}-1,\frac{r_n}{r_{n+1}}\Big].$$
So we can partition every square of side length $r_n$ into $r_n^2/r_{n+1}^2$ squares of side length $r_{n+1}$.
Using property (e) of $x$, it follows easily by induction that
$$\frac{1}{\sqrt{2}(E^n)'(x)}<r_n\text{ for all }n\in\mathbb{N}.$$
Consequently,
\begin{equation}\label{84}
\frac{1}{\sqrt{2}(E^n)'(x)}<r_n\leq\frac{1}{(E^n)'(x)}\text{
for all } n\in\mathbb{N}.
\end{equation}
Choose $n\in\mathbb{N}$ and a square $Q=Q(z_0,r_n,0)\subset
Q_0$ such that $Q\cap I_R(E)\neq\emptyset$.

\emph{Claim}: If we partition $Q$ into squares of side length $r_{n+1}$, then at least $$\Big(\frac{r_n}{r_{n+1}}\Big)^2\Big(\frac{1}{2}-\eps'\Big)$$ of these squares are contained in $\mathcal{F}(E)$.

We start the proof of this claim by showing that there exists $k\in\mathbb{N}$ such that $E^k(Q)$ is 'close to a square' of side length bigger than $Mr_0$.

Because of \eqref{78}, we have dist$(Q,P(E))\geq C$. By
the blow-up property of $\mathcal{J}(E)$ (Lemma \ref{blow}), there exists a minimal $k\in\mathbb{N}$
such that $E^{k-1}(Q)$ is contained in some square
$Q(E^{k-1}(z_0),Mr_0)$, but $E^k(Q)$ is not. We
now show that $E^{k-1}:Q\to\Com$ has a univalent
continuation to $D(z_0,Kr_n/\sqrt{2})$. Let $\varphi$ be the
branch of $(E^{k-1})^{-1}$ which is defined on a
neighborhood of $E^{k-1}(z_0)$ and satisfies
$\varphi(E^{k-1}(z_0))=z_0$. By the monodromy theorem and
since
$$\text{dist}(E^{k-1}(z_0),P(E))\geq\text{dist}(\set{\Re z\geq R},P(E))-\sqrt{2}Mr_0\geq
C>4KMr_0,$$ we conclude from \eqref{78} that $\varphi$ has an
analytic continuation to $D(E^{k-1}(z_0),4KMr_0)$. Because
$E^{k-1}(Q)$ is contained in a disc of radius
$Mr_0/\sqrt{2}$, it follows from Lemma \ref{cif} that
$$\abs{(E^{k-1})'(z_0)}\leq\frac{Mr_0/\sqrt{2}}{r_n/2}=M\sqrt{2}\frac{r_0}{r_n},$$
and consequently
$$\abs{\varphi'(E^{k-1}(z_0))}\geq\frac{1}{Mr_0}\frac{r_n}{\sqrt{2}}.$$
By the Koebe 1/4-Theorem, we obtain
$$\varphi(D(E^{k-1}(z_0),4KMr_0))\supset
D(z_0,\abs{\varphi'(E^{k-1}(z_0))}KMr_0)\supset
D(z_0,Kr_n/\sqrt{2}),$$ and therefore $E^{k-1}$ can be
continued univalently to $D(z_0,Kr_n/\sqrt{2})$. By Lemma
\ref{first} and \eqref{103}, the distortion of $E^{k-1}$
satisfies
\begin{equation}\label{80}
D(E^{k-1}|_Q)\leq\Big(\frac{K+1}{K-3}\Big)^6<1+2\eps.
\end{equation}
It follows from Lemma \ref{second} that $E^{k-1}(Q)$
contains some square
\begin{align*}
Q_1&:=Q(E^{k-1}(z_0),\abs{(E^{k-1})'(z_0)}r_n(1-\sqrt{2}\eps)/(1+2\eps))\\
&\supset
Q(E^{k-1}(z_0),\abs{(E^{k-1})'(z_0)}r_n(1-2\eps)^2)
\end{align*} and is contained in some square
\begin{align*}
Q_2&:=Q(E^{k-1}(z_0),\abs{(E^{k-1})'(z_0)}r_n(1+\sqrt{2}\eps)(1+2\eps))\\
&\subset
Q(E^{k-1}(z_0),\abs{(E^{k-1})'(z_0)}r_n(1+2\eps)^2).
\end{align*}
Because the side length $l_1$ of $Q_1$ is at most $Mr_0$, the side
length $l_2$ of $Q_2$ satisfies
$$l_2<\Big(\frac{1+2\eps}{1-2\eps}\Big)^2Mr_0.$$
Since
$$\Big(\frac{1+2\eps}{1-2\eps}\Big)^2KMr_0\frac{1}{\sqrt{2}}<\pi$$
by \eqref{79}, $E$ is univalent on
$$D\Big(E^{k-1}(z_0),\Big(\frac{1+2\eps}{1-2\eps}\Big)^2KMr_0/\sqrt{2}\Big)\supset
D\Big(E^{k-1}(z_0),\frac{l_2K}{\sqrt{2}}\Big).$$ The same arguments as before yield
\begin{equation}\label{81}
D(E|_{Q_2})<1+2\eps,
\end{equation}
and hence $E(Q_1)$ contains a square of side length
$\abs{E'(E^{k-1}(z_0))}l_1(1-\sqrt{2}\eps)/(1+2\eps)$
and $E(Q_2)$ is contained in a square of side length
$\abs{E'(E^{k-1}(z_0))}l_2(1-\sqrt{2}\eps)(1+2\eps)$.
Thus, summarizing the previous estimates, $E(Q_1)$
contains some square
$$Q_1':=Q\Big(E^k(z_0),\abs{(E^k)'(z_0)}r_n\frac{(1-\sqrt{2}\eps)(1-2\eps)^2}{1+2\eps}\Big)$$
and $E(Q_2)$ is contained in some square
$$Q_2':=Q\Big(E^k(z_0),\abs{(E^k)'(z_0)}r_n(1+\sqrt{2}\eps)(1+2\eps)^3\Big).$$
In particular, $Q_1'\subset E^k(Q)\subset Q_2'$.
Since $E^k(Q)$ is not contained in a square of side length
$Mr_0$, it follows that the side length $l_2'$ of $Q_2'$ satisfies
$$l_2'>Mr_0.$$
By the definition of $a$ (see \eqref{1234}) and \eqref{12345}, we conclude that
$$l_1'=\frac{l_2'}{1+a}>\frac{1}{1+a}Mr_0\geq(M-1)r_0.$$
If $\eps$ is chosen small enough, we can assume that the side
length $l_1'$ of $Q_1'$ is also at least $Mr_0$. There are two
cases:

\emph{Case 1}: $l_1'\leq Ar_0$. \\ Then (if $A$ is chosen large
enough) it follows that $l_2'\leq(A+1)r_0$.

We now show that at least almost half of $Q$ belongs to the Fatou set.
Choose $m\in\mathbb{N}$ such that $mr_0\leq r_1'$, but $(m+1)r_0>r_1'$.
Then $m\geq M-1$ and we can find $m^2r_0^2$ squares of side length $r_0$
which have pairwise disjoint interiors and which are
completely contained in $Q_1'$. Let $P=Q(p,r_0)$ be one of these
squares. Because $\Re z\geq R-(A+1)r_0\sqrt{2}\geq R_1$ for
all $z\in P$ and by the way $R_1$ and $r_0$ are defined, there
exists an $r_0$-packing pack$_{E(P)}$ of $E(P)\cap
F_\delta$ such that
$$\text{dens}(\text{pack}_{E(P)},E(P))>\frac{1}{2}-\eps.$$
Since $E$ is clearly injective on $D(p,KMr_0/\sqrt{2})$
(because $KMr_0/\sqrt{2}<\pi$), the distortion satisfies
$$D(E|_P)<1+2\eps.$$ If $\varphi_P$ denotes the branch of
$E^{-1}$ such that
$\varphi_P(\text{pack}_{E(P)})\subset P$, we obtain
by Lemma \ref{eigensch} (d) that
$$\text{dens}\big(\varphi_P(\text{pack}_{E(P)}),P\big)\geq\frac{1}{(1+2\eps)^2}\Big(\frac{1}{2}-\eps\Big).$$
If we take the union over all $P$, it follows that
$$\text{dens}\Big(\bigcup_P\varphi_P(\text{pack}_{E(P)}),Q_1'\Big)\geq\frac{1}{(1+2\eps)^2}\Big(\frac{1}{2}-\eps\Big)\Big(\frac{m}{m+1}\Big)^2$$
and consequently
\begin{align*}
\text{dens}\Big(\bigcup_P\varphi_P(\text{pack}_{E(P)}),E^k(Q)\Big)&\geq\frac{1}{(1+2\eps)^2}\Big(\frac{1}{2}-\eps\Big)\Big(\frac{m}{m+1}\Big)^2\frac{\abs{Q_1'}}{\abs{Q_2'}}\\
&=\Big(\frac{m}{m+1}\Big)^2\Big(\frac{1}{2}-\eps\Big)\frac{(1-2\eps)^4(1-\sqrt{2}\eps)^2}{(1+2\eps)^{10}(1+\sqrt{2}\eps)^2}.
\end{align*}
Taking preimages, it follows from \eqref{80} and \eqref{81} that
\begin{align}\label{96}
&&&\hspace{0.2cm}\text{dens}\Big((E^k)^{-1}\Big(\bigcup_P\varphi_P(\text{pack}_{E(P)})\Big)\cap
Q,Q\Big)\notag\\&&\geq&\hspace{0.2cm}\frac{1}{D(E^k|_Q)^2}\Big(\frac{m}{m+1}\Big)^2\Big(\frac{1}{2}-\eps\Big)\frac{(1-2\eps)^4(1-\sqrt{2}\eps)^2}{(1+2\eps)^{10}(1+\sqrt{2}\eps)^2}\notag\\
&&\geq&\hspace{0.2cm}\frac{1}{D(E^{k-1}|_Q)^2}\frac{1}{D(E|_{Q_2})^2}\Big(\frac{m}{m+1}\Big)^2\Big(\frac{1}{2}-\eps\Big)\frac{(1-2\eps)^4(1-\sqrt{2}\eps)^2}{(1+2\eps)^{10}(1+\sqrt{2}\eps)^2}\notag\\
&&\geq&\hspace{0.2cm}\Big(\frac{M-1}{M}\Big)^2\Big(\frac{1}{2}-\eps\Big)\frac{(1-2\eps)^4(1-\sqrt{2}\eps)^2}{(1+2\eps)^{14}(1+\sqrt{2}\eps)^2}
\end{align}
We show now that these preimages (which all belong to the Fatou
set of $E$) are 'large' compared to a square of side
length $r_{n+1}$.

Let $P$ be a square in $Q_1'$, and let $z\in P$
be arbitrary. Then $\Re z\geq R-(A+1)\sqrt{2}r_0$, and consequently
$$\abs{E(z)}\geq\lambda'e^{R-(A+1)\sqrt{2}r_0},$$
which implies that
$$\text{dist}(E(z),P(E))\geq\abs{E(z)}-\alpha>C$$
by the definition of $R$. Hence the branch $\psi$ of
$(E^{k+1})^{-1}$ with $\psi(E(z))\in Q$ can be
continued univalently to $D(E(z),Kr_0/\sqrt{2})$. We
denote (as in section 3) the collection of boxes that form
pack$_{E(P)}$ by $\mathcal{B}(E(P))$. So if
$P'=Q(z_{P'},r_0,0)\in\mathcal{B}(E(P))$ and $\psi_{P'}$
denotes the branch of $(E^{k+1})^{-1}$ with
$\psi_{P'}(P')\subset Q$, then $D(\psi_{P'}|_{P'})<1+2\eps$ and
$\psi_{P'}(P')$ is contained in a square $Q_{P'}'$ of side length
\begin{equation}\label{def1}
\rho_{P'}':=\abs{\psi_{P'}'(z_{P'})}r_0(1+2\eps)(1+\sqrt{2}\eps)
\end{equation}
and
contains a square $Q_{P'}$ of side length
\begin{equation}\label{def2}
\rho_{P'}:=\abs{\psi_{P'}'(z_{P'})}r_0\frac{1}{1+2\eps}(1-\sqrt{2}\eps).
\end{equation}
Let
\begin{equation}\label{def3}
w_{P'}:=\psi_{P'}(z_{P'})\in Q.
\end{equation}
We now show that
\begin{equation}\label{85}
(E^{n+1})'(x)\geq L\abs{(E^{k+1})'(w_{P'})},
\end{equation}
which means that $r_{n+1}$ is much smaller than the size of the
preimages of pack$_{E(P)}$ in $Q$. If $k\leq n+2$, then
\eqref{85} is clear by \eqref{87}, so we may assume that $k-n\geq
3.$ Since $C>90(A+1)r_0$, we obtain
$$\text{dist}(E^k(Q),P(E))\geq\text{dist}(Q_2',P(E))>6(A+1)r_0.$$
Hence the branch $\varphi_0$ of $(E^k)^{-1}$ with
$\varphi_0(E^k(z_0))=z_0$ can be continued analytically to
$D(E^k(z_0),6(A+1)r_0)$. Because
\begin{equation}\label{83}
\abs{(E^k)'(z_0)}\leq\sqrt{2}(A+1)\frac{r_0}{r_n}
\end{equation}
by the same arguments as above, it follows again by the Koebe
1/4-Theorem that
$$\varphi_0(D(E^k(z_0),6(A+1)r_0))\supset
D\Big(z_0,\frac{3}{2}\frac{r_n}{\sqrt{2}}\Big).$$ Thus
$E^k$ is univalent on $D(z_0,\frac{3r_n}{2\sqrt{2}})$.
Corollary \ref{kk} implies that
\begin{equation}\label{82}
\abs{(E^k)'(w_{P'})}\leq\frac{3r_n/2\sqrt{2}+r_n/\sqrt{2}}{(3r_n/2\sqrt{2}-r_n/\sqrt{2})^3}
(3r_n/2\sqrt{2})^2\abs{(E^k)'(z_0)}.
\end{equation}
Combining \eqref{84}, \eqref{83} and \eqref{82}, we obtain
\begin{equation}\label{92}
\abs{E^k(w_{P'})}\leq\frac{5/2(3/2)^2}{(1/2)^3}\frac{\sqrt{2}(A+1)r_0}{r_n}=\frac{B}{\sqrt{2}r_n}\leq
B\, (E^n)'(x).
\end{equation}
For the proof of \eqref{85}, it remains to show that
\begin{equation}\label{86}
\kappa E^{n+1}(x)\geq\abs{E^{k+1}(w_{P'})}.
\end{equation}
To see this, we first show that
\begin{equation}\label{89}
x^{\kappa^{1/2^j}-1}\leq\kappa^{1/2^{j-1}}\text{ for all
}j\in\mathbb{N}_0,
\end{equation}
which is equivalent to
$$x\geq\frac{x^{(\kappa^{1/2^j})}}{\kappa^{1/2^{j-1}}}\text{ for all }j\in\mathbb{N}_0.$$
Because of the obvious inequality $$2^j(\kappa^{1/2^{j+1}}-\kappa^{1/2^j})\geq\sqrt{\kappa}-\kappa\text{
for all }j\in\mathbb{N}_0$$ and condition (b) on $x$, the
right side is an increasing sequence tending to $x$ as $j\to\infty$, and \eqref{89} is proven.
The inequality \eqref{89} easily implies
\begin{equation}\label{90}
E^{i}(x)^{\kappa^{1/2^j}-1}\leq\kappa^{1/2^{j-1}} \text{
for all }i,j\in\mathbb{N}_0.
\end{equation}
Now suppose that
\begin{equation}\label{91}
\kappa E^{n+1}(x)<\abs{E^{k+1}(w_{P'})}.
\end{equation}
We show that this implies
$$\kappa^{1/2^j}E^{n+1-j}(x)<\abs{E^{k+1-j}(w_{P'})}\text{ for all }1\leq j\leq n.$$
This can be done inductively. For $j=1$, suppose that
$\abs{E^k(w_{P'})}\leq\kappa^{1/2}E^n(x).$ Then
$$\abs{E^{k+1}(w_{P'})}=\lambda'e^{\Re E^k(w_{P'})}
\leq\lambda'e^{\kappa^{1/2}E^n(x)}
\leq(\lambda'e^{E^n(x)})^{\kappa^{1/2}}
=E^{n+1}(x)^{\kappa^{1/2}}
\leq\kappa E^{n+1}(x),$$
where the second inequality holds because $\lambda'<1$, and the
last inequality follows from \eqref{90}. The step $j\to j+1$
can be done in exactly the same way. So \eqref{91} implies
$$\kappa^{1/2+\ldots+1/2^n}\prod_{j=1}^nE^j(x)<\hspace{-0.4cm}\prod_{j=k-n+1}^k\hspace{-0.4cm}\abs{E^j(w_{P'})},$$
which is the same as
\begin{equation}\label{95}
\kappa^{1-1/2^n}(E^n)'(x)<\hspace{-0.4cm}\prod_{j=k-n+1}^k\hspace{-0.4cm}\abs{E^j(w_{P'})}.
\end{equation}
We have
\begin{align}\label{93}
\prod_{j=k-n+1}^k\hspace{-0.4cm}\abs{E^j(w_{P'})}&=\frac{1}{\prod_{j=1}^{k-n}\abs{E^j(w_{P'})}}\abs{(E^k)'(w_{P'})}\notag\\
&\leq\frac{B}{\prod_{j=1}^{k-n}\abs{E^j(w_{P'})}}(E^n)'(x)
\end{align}
by \eqref{92}, and
\begin{equation}\label{94}
\abs{E^j(w_{P'})}>R-Mr_0\sqrt{2}>C-Mr_0\sqrt{2}=L\text{
for }j\leq k-n
\end{equation}
since $I_R(E)\cap Q\neq\emptyset$ and $E^j(Q)$ is contained in a
square of side length $Mr_0$ for $j\leq k-1$.
Combining \eqref{95}, \eqref{93} and \eqref{94} yields
$$\kappa^{1-1/2^n}(E^n)'(x)<\frac{B}{L^{k-n}}(E^n)'(x),$$
and hence
$$\frac{L^{k-n}\kappa^{1-1/2^n}}{B}<1.$$ We now can show easily
that this is impossible: It suffices to prove that
$$\frac{L^{k-n}\kappa}{B}=\frac{L^{k-n-1}}{B^2}\geq 1.$$
But this is clear since $k-n-1\geq 2$ and $L>B$ by our assumption.
So we have finally proved \eqref{85}. Consequently, the side length $\rho_{P'}$ of the square $Q_{P'}$ satisfies $$\rho_{P'}\geq
L\frac{1}{(E^{n+1})'(x)}r_0\frac{1}{1+2\eps}(1-\sqrt{2}\eps)=N\frac{1}{(E^{n+1})'(x)}\geq
Nr_{n+1}.$$
Now we are ready to prove the claim in this case.
At least $\big(\frac{\rho_{P'}}{r_{n+1}}\big)^2$
squares in the $r_{n+1}$-grid intersect $Q_{P'}$, whereas only
$C'\frac{\rho_{P'}}{r_{n+1}}$ squares of the grid intersect
$\partial Q_{P'}$, which implies that
$$\Big(\frac{\rho_{P'}}{r_{n+1}}\Big)^2-C'\frac{\rho_{P'}}{r_{n+1}}=\Big(\frac{\rho_{P'}}{r_{n+1}}\Big)^2\Big(1-C'\frac{r_{n+1}}{\rho_{P'}}\Big)
\geq\Big(\frac{\rho_{P'}}{r_{n+1}}\Big)^2\Big(1-C'\frac{1}{N}\Big)>\Big(\frac{\rho_{P'}}{r_{n+1}}\Big)^2(1-\eps)$$
squares are contained in $Q_{P'}$. Because all the previous
estimates were independent of $P'$ and $P$, we can do the same for
every $P$ and every $P'\in\mathcal{B}(E(P))$ and obtain
that at least
$$(1-\eps)\sum_P\sum_{P'\in\mathcal{B}(E(P))}\Big(\frac{\rho_{P'}}{r_{n+1}}\Big)^2$$
squares of the $r_{n+1}$-grid are contained in $\mathcal{F}(E)$. It follows now from \eqref{96} that
\begin{align*}
&&&\hspace{0.2cm}(1-\eps)\sum_P\sum_{P'\in\mathcal{B}(E(P))}\Big(\frac{\rho_{P'}}{r_{n+1}}\Big)^2\\
&&=&\hspace{0.2cm}\frac{1}{r_{n+1}^2}(1-\eps)\frac{1}{(1+2\eps)^4}\Big(\frac{1-\sqrt{2}\eps}{1+\sqrt{2}\eps}\Big)^2\sum_P\sum_{P'\in\mathcal{B}(E(P))}\hspace{-0.4cm}(\rho_{P'}')^2\\
&&\geq&\hspace{0.2cm}\frac{1}{r_{n+1}^2}(1-\eps)\frac{1}{(1+2\eps)^4}\Big(\frac{1-\sqrt{2}\eps}{1+\sqrt{2}\eps}\Big)^2\sum_P\sum_{P'\in\mathcal{B}(E(P))}\hspace{-0.4cm}\abs{\psi_{P'}(P')}^2\\
&&=&\hspace{0.2cm}\Big(\frac{r_n}{r_{n+1}}\Big)^2(1-\eps)\frac{(1-2\eps)^4}{(1+2\eps)^{18}}\Big(\frac{1-\sqrt{2}\eps}{1+\sqrt{2}\eps}\Big)^4\Big(\frac{M-1}{M}\Big)^2\Big(\frac{1}{2}-\eps\Big)\\
&&\geq&\hspace{0.2cm}\Big(\frac{r_n}{r_{n+1}}\Big)^2\Big(\frac{1}{2}-\eps'\Big),
\end{align*}
so at least
$$\Big(\frac{r_n}{r_{n+1}}\Big)^2\Big(\frac{1}{2}-\eps'\Big)$$
squares from the $r_{n+1}$-grid belong to
$\mathcal{F}(E)$.

\emph{Case 2:} $l_1'>Ar_0$.
In this case, we can find an $r_0$-packing pack$_{Q_1'}$ of
$Q_1'\cap F_\delta$ such that
$$\text{dens}(\text{pack}_{Q_1'},Q_1')>\frac{1}{2}-\eps.$$
Because we have to be sure that $(E^k)^{-1}$ is univalent
on a domain which is large enough, we consider only the squares of
pack$_{Q_1'}$ which are contained in $Q_1'\setminus Y$, where
$Y:=\set{z\in\Com\trenn\text{dist}(z,P(E))<\pi}\subset
D(0,\pi+1)$. By the definition of $A$, we have
$$\text{dens}(Y,Q_1')<\eps,$$
hence if pack$_{Q_1'\setminus Y}$ denotes only the elements of
pack$_{Q_1'}$ which are contained in $Q_1'\setminus Y$, we easily
obtain
$$\text{dens}(\text{pack}_{Q_1'\setminus Y},Q_1')\geq\text{dens}(\text{pack}_{Q_1'},Q_1')-\text{dens}(Y,Q_1')>\frac{1}{2}-2\eps.$$
Because $D(E^k|_Q)<(1+2\eps)^2$, it follows that
\begin{align}\label{98}
\text{dens}((E^k)^{-1}(\text{pack}_{Q_1'\setminus Y})\cap
Q,Q)&\geq\frac{1}{(1+2\eps)^4}\text{dens}(\text{pack}_{Q_1'\setminus Y},E^k(Q))\notag\\
&=\frac{1}{(1+2\eps)^4}\text{dens}(\text{pack}_{Q_1'\setminus Y},Q_1')\text{dens}(Q_1',E^k(Q))\notag\\
&\geq\frac{1}{(1+2\eps)^4}\Big(\frac{1}{2}-2\eps\Big)\text{dens}(Q_1',Q_2')\notag\\
&=\frac{(1-2\eps)^4}{(1+2\eps)^{12}}\frac{(1-\sqrt{2}\eps)^2}{(1+\sqrt{2}\eps)^2}\Big(\frac{1}{2}-2\eps\Big).
\end{align}
Again, denote by $\mathcal{B}(Q_1'\setminus Y)$ the collection of
$r_0$-boxes that form pack$_{Q_1'\setminus Y}$ and let
$P'=Q(z_{P'},r_0,0)\in\mathcal{B}(Q_1'\setminus Y)$. By the
definition of $Y$, it follows that the branch $\psi_{P'}$ of
$(E^k)^{-1}$ that maps $P'$ into $Q$ can be continued
univalently to $D(z_R,KMr_0/\sqrt{2})$. Hence we have as before
that $\psi(P')$ is contained in a square of side length $\rho'_{P'}$ and
contains a square $Q_{P'}$ of side length $\rho_{P'}$, where $\rho'_{P'}$ and $\rho_{P'}$
are defined as in \eqref{def1} and \eqref{def2}. Let
$w_{P'}:=\psi_{P'}(z_{P'})\in Q$ be as in \eqref{def3}. We now show that
\begin{equation}\label{97}
(E^{n+1})'(x)\geq L\abs{(E^k)'(w_{P'})}.
\end{equation}
First note that because $E^{k-1}$ is injective on
$D(z_0,Kr_n/\sqrt{2})$, we conclude as in case~1 (using Corollary
\ref{kk} again) that
$$\abs{(E^{k-1})'(w_{P'})}\leq B\, (E^n)'(x),$$ and by
the same arguments as before, \eqref{97} follows. Again we see
(analogously to case 1) that $Q_{P'}$ contains at least
$(1-\eps)\big(\frac{\rho_{P'}}{r_{n+1}}\big)^2$ squares in the
$r_{n+1}$-grid. Since all of our estimates were
independent of $P'$, we obtain that at least
$$(1-\eps)\hspace{-0.3cm}\sum_{P'\in\mathcal{B}(Q_1'\setminus Y)}\hspace{-0.3cm}\Big(\frac{\rho_{P'}}{r_{n+1}}\Big)^2$$
squares in the $r_{n+1}$-grid are contained in
$\mathcal{F}(E)$, and thus we have with \eqref{98} that
\begin{align*}
(1-\eps)\sum_{P'\in\mathcal{B}(Q_1'\setminus Y)}\Big(\frac{\rho_{P'}}{r_{n+1}}\Big)^2&=\frac{1}{r_{n+1}^2}(1-\eps)\hspace{-0.4cm}\sum_{P'\in\mathcal{B}(Q_1'\setminus Y)}\hspace{-0.4cm}(\rho_{P'})^2\\
&=\frac{1}{r_{n+1}^2}(1-\eps)\hspace{-0.4cm}\sum_{P'\in\mathcal{B}(Q_1'\setminus Y)}\hspace{-0.4cm}(\rho_{P'}')^2\Big(\frac{1-\sqrt{2}\eps}{1+\sqrt{2}\eps}\Big)^2\frac{1}{(1+2\eps)^4}\\
&\geq\frac{1}{r_{n+1}^2}(1-\eps)\Big(\frac{1-\sqrt{2}\eps}{1+\sqrt{2}\eps}\Big)^2\frac{1}{(1+2\eps)^4}\hspace{-0.2cm}\sum_{P'\in\mathcal{B}(Q_1'\setminus Y)}\hspace{-0.4cm}\abs{\psi_{P'}(P')}\\
&\geq\frac{1}{r_{n+1}^2}(1-\eps)\Big(\frac{1-\sqrt{2}\eps}{1+\sqrt{2}\eps}\Big)^4\frac{(1-2\eps)^4}{(1+2\eps)^{16}}\Big(\frac{1}{2}-2\eps\Big)r_n^2\\
&\geq\Big(\frac{r_n}{r_{n+1}}\Big)^2\Big(\frac{1}{2}-\eps'\Big),
\end{align*}
and our claim is also true in case 2.

So in each of the two cases,
$$\Big(\frac{1}{2}+\eps'\Big)\Big(\frac{r_n}{r_{n+1}}\Big)^2$$ squares of side length
$r_{n+1}$ suffice to cover $I_R(E)\cap Q$. Hence we obtain the following: Starting with a
square $Q$ of side length $r_1$ depending on $Q$, we can cover $I_R(E)\cap Q$ by
$$\Big(\frac{1}{2}+\eps'\Big)^n\Big(\frac{r_1}{r_{n+1}}\Big)^2$$
squares of side length $r_{n+1}$. It follows (with
$h_\gamma(t)=t^2\Phi(1/t)^\gamma$) that
\begin{align*}
\mathcal{H}^{h_\gamma}(I_R(E)\cap
Q)&\leq\lim_{n\to\infty}\Big(\frac{1}{2}+\eps'\Big)^n\Big(\frac{r_1}{r_{n+1}}\Big)^2h_\gamma(\sqrt{2}r_{n+1})\\
&=\lim_{n\to\infty}\Big(\frac{1}{2}+\eps'\Big)^nr_1^22\Phi\Big(\frac{1}{\sqrt{2}r_{n+1}}\Big)^\gamma\\
&<2r_1^2\lim_{n\to\infty}\Big(\frac{1}{2}+\eps'\Big)^n\Phi((E^{n+1})'(x))^\gamma\\
&\leq2r_1^2\lim_{n\to\infty}\Big(\frac{1}{2}+\eps'\Big)^n\Phi(E^{n+1}(2x))^\gamma\\
&=2r_1^2\lim_{n\to\infty}\Big(\frac{1}{2}+\eps'\Big)^n(\beta^{n+1})^\gamma\Phi(2x)^\gamma\\
&=2r_1^2\beta^\gamma\Phi(2x)^\gamma\lim_{n\to\infty}\Big(\Big(\frac{1}{2}+\eps'\Big)\beta^\gamma\Big)^n\\
&=0,
\end{align*}
because $\beta^\gamma(1/2+\eps')<1$. We now only need to
show that $I_R(E)$ can be covered by countably many
such squares $Q$: For $Q_0=Q(z_0,r_0,0)$,
choose $x=x(Q_0)$ as before (which defines $r_1$), and the intersection of $I_R(E)$ with every square
of side length $r_1$ which is completely contained in $Q_0$ has
zero $\mathcal{H}^{h_\gamma}$-measure. Because $Q_0$ can
be covered by finitely many squares of side length $r_1$ and
$I_R(E)$ can clearly be covered by countably many squares
of side length $r_0$, the theorem is proved.
\end{proo2}

Observe that because $I(E)\subset\mathcal{J}(E)$, the previous
theorem is clearly also true when we replace
$\mathcal{J}(E)$ by $I(E)$.

\section{The main theorem and its consequences}

\subsection{Equivalence of the gauge functions}

In the preceding two sections, we estimated
$\mathcal{H}^{h_\gamma}(\mathcal{J}(E_\lambda))$ for $\gamma>0$, where
$h_\gamma=h_{\lambda',\gamma}$ with fixed $\lambda'\in(0,1/e)$.
One major difference between section 3 and section 4 is that $\lambda$ was an
arbitrary element of $\Com\setminus\set{0}$ in section 3, whereas $\lambda$ had to equal
$\lambda'$ in section 4. In this section, we will
prove that if $\lambda\in(0,1/e)$, then $\mathcal{J}(E_\lambda)$
still has zero measure with respect to the gauge function
$h_\gamma$, where $\gamma$ is chosen such that $\beta^\gamma<2$.
We prove this by showing that two gauge functions
$h_{\lambda_1,\gamma_1}$ and $h_{\lambda_2,\gamma_2}$ are
equivalent (in the sense that there exist constants $c,C>0$ with
$c\, h_{\lambda_1,\gamma_1}(t)\leq
h_{\lambda_2,\gamma_2}(t)\leq C\, h_{\lambda_1,\gamma_1}(t)$
if $t$ is small) whenever
$\beta_{\lambda_1}^{\gamma_1}=\beta_{\lambda_2}^{\gamma_2}$.

To prove this result, it is more convenient to consider a
different parametrization of the exponentials with a real repelling fixed point. Let
$\lambda\in(0,1/e)$ be fixed. If $\mu\in(1,\infty)$ is chosen
such that $\mu/e^\mu=\lambda$, i.e.
$\mu=\beta_\lambda$, we conjugate $E_\lambda$ by $z\mapsto
z+\mu$ and get a corresponding function
$$\tilde{E}_\mu(z)=\mu(e^z-1).$$ Then 0 is a repelling fixed
point of $\tilde{E}_\mu$ with multiplier $\mu$. If we denote the
Poincaré function of $\tilde{E}_\mu$ with respect to 0 by $\tilde{\Phi}_\mu$, it is easy
to check that $\tilde{\Phi}_\mu(z)=\Phi_\lambda(z+\mu)$.
If we further set $\tilde{h}_{\mu,\gamma}(t):=t^2(\tilde{\Phi}_\mu(1/t))^\gamma$, then
it is not hard to see that the gauge functions $\tilde{h}_{\mu,\gamma}$
and $h_{\lambda,\gamma}$ define the same Hausdorff measure.

Consequently, the desired statement follows from
\begin{theo}
Let $1<\mu_1,\mu_2<\infty$ and $\gamma_1,\gamma_2>0$ such that
$\mu_1^{\gamma_1}=\mu_2^{\gamma_2}$. Let
$\tilde{\Phi}_{\mu_1},\tilde{\Phi}_{\mu_2}$ and
$\tilde{h}_{\mu_1,\gamma_1},\tilde{h}_{\mu_2,\gamma_2}$ be
defined as above. Then there exist constants $c,C>0$ such that
$$c<\frac{\tilde{h}_{\mu_1,\gamma_1}(t)}{\tilde{h}_{\mu_2,\gamma_2}(t)}<C$$
for every $t>0$ small enough.
\end{theo}
\begin{proo}
First, for an arbitrary $\mu\in(1,\infty)$, define
$$L_\mu:=(\tilde{E}_\mu|_\mathbb{R})^{-1}.$$ Then $L_\mu$ is
defined on $(-\mu,\infty)$ and a simple computation shows that
\begin{equation*}
L_\mu(x)=\log\Big(1+\frac{x}{\mu}\Big).
\end{equation*}
Choose $x_1(\mu)$ with
\begin{equation}\label{71}
L_\mu(x)\leq\log x\text{ whenever }x\geq x_1(\mu)
\end{equation}
and
\begin{equation*}
x_0(\mu):=L_\mu(x_1(\mu))>1.
\end{equation*}
The functional equation for $\tilde{\Phi}_\mu$ gives
\begin{equation*}
\tilde{\Phi}_\mu(L_\mu^n(x))=\frac{1}{\mu^n}\tilde{\Phi}_\mu(x)\text{
for all }n\in\mathbb{N}.
\end{equation*}
Let $x\geq x_1(\mu)$ and $n_\mu(x)\in\mathbb{N}$ such that
\begin{equation*}
L_\mu^{n_\mu(x)}(x)\in[x_0(\mu),x_1(\mu)).
\end{equation*}
Then
\begin{align}\label{77}
\tilde{\Phi}_\mu(x)&=\lim_{n\to\infty}\mu^nL_\mu^n(x)\notag\\
&=\lim_{n\to\infty}\mu^nL_\mu^{n-n_\mu(x)}(L_\mu^{n_\mu(x)}(x))\notag\\
&=\mu^{n_\mu(x)}\lim_{n\to\infty}\mu^{n-n_\mu(x)}L_\mu^{n-n_\mu(x)}(L_\mu^{n_\mu(x)}(x))\notag\\
&=\mu^{n_\mu(x)}\tilde{\Phi}_\mu(L_\mu^{n_\mu(x)}(x)).
\end{align}
By induction, it follows easily from \eqref{71} and the monoticity
of the logarithm that
\begin{equation}\label{72}
L_\mu^n(x)\leq\log^n(x)\text{ if }n\leq n_\mu(x).
\end{equation}
We now prove (also by induction) an inequality in the other
direction: We claim that
\begin{equation}\label{73}
L_\mu^n(x)\geq\log^n(x)-\log(\mu)\text{ if }n\leq n_\mu(x).
\end{equation}
First of all, observe that
$$L_\mu(x)=\log\Big(1+\frac{x}{\mu}\Big)=\log\Big(\frac{x}{\mu}\Big)+\log\Big(1+\frac{\mu}{x}\Big)>\log x-\log\mu,$$
so the inequality is true for $n=1$. Now suppose that it holds for
some $n<n_\mu(x)$. Then it follows from \eqref{72} that
$\log^{n+1}(x)\geq x_0(\mu)>1$, in particular $\log^{n+1}(x)$ is
well defined. We compute
\begin{align*}
L_\mu^{n+1}(x)&=L_\mu(L_\mu^n(x))\\
&>L_\mu(\log^n(x)-\log(\mu))\\
&=\log\Big(1+\frac{\log^n(x)-\log(\mu)}{\mu}\Big)\\
&=\log^{n+1}(x)-\log(\mu)+\log\Big(1+\frac{\mu}{\log^n(x)}-\frac{\log(\mu)}{\log^n(x)}\Big)\\
&>\log^{n+1}(x)-\log(\mu),
\end{align*}
which finishes the proof of \eqref{73}. We have thus proved that
for every $\mu\in(1,\infty)$,
\begin{equation}\label{74}
\log^n(x)-\log\mu\leq L_\mu^n(x)\leq\log^n(x)
\end{equation}
for all $x\geq x_1(\mu)$ and $n\leq n_\mu(x)$. \\
Let such $x_1(\mu_1)$ and $x_1(\mu_2)$ be chosen for $\mu_1$
and $\mu_2$ and define
$$x_1:=\max\set{x_1(\mu_1),x_1(\mu_2)}.$$ We show that the
difference between $n_{\mu_1}(x)$ and $n_{\mu_2}(x)$ is
uniformly bounded in $x$. Let $x\geq x_1$. Without loss of
generality we may assume that $n_{\mu_1}(x)\geq n_{\mu_2}(x)$.
Using \eqref{74}, we obtain
$$x_0(\mu_1)\leq\log^{n_{\mu_1}(x)}(x)\leq
L_\mu^{n_{\mu_1}(x)}(x)+\log\mu_1\leq x_1+\log\mu_1$$ as
well as the corresponding statement for $\mu_2$. If we define
$x_0:=\min\set{x_0(\mu_1),x_0(\mu_2)}$, then
$$\log^{n_{\mu_j}(x)}(x)\in[x_0,x_1+\log\max_i\mu_i]\text{ for all }x\geq x_1,j=1,2.$$
Now suppose that $(x_k)$ is a sequence of real numbers tending to
$\infty$ such that
$$\abs{(n_{\mu_1}-n_{\mu_2})(x_k)}\to\infty.$$ Because
$\log^{n_{\mu_1}(x_k)}(x_k)$ is bounded below by $x_0$ and
$\log^{n_{\mu_2}(x_k)}(x_k)$ is bounded above by
$M=x_1+(1+\eps)\log\max_i\mu_i$, we obtain
$$x_0\leq\log^{n_{\mu_1}(x_k)}(x_k)=\log^{(n_{\mu_1}-n_{\mu_2})(x_k)}(\log^{n_{\mu_2}(x_k)}(x_k))\leq\log^{(n_{\mu_1}-n_{\mu_2})(x_k)}(M),$$
which is clearly impossible. It follows that there exists some
constant $K$ with
\begin{equation}\label{76}
\abs{n_{\mu_1}(x)-n_{\mu_2}(x)}\leq K\text{ for all }x\geq x_1,
\end{equation}
thus this difference is uniformly bounded in $x$. Let
$n(x):=n_{\mu_1}(x).$ We use \eqref{77} to deduce that for
$i=1,2$ and $x\geq x_1$,
$$\tilde{\Phi}_{\mu_i}(x)
\begin{cases}
\leq\mu_i^{n_{\mu_i}(x)}\tilde{\Phi}_{\mu_i}(x_1)=:\mu_i^{n_{\mu_i}(x)}B_i\\
\geq\mu_i^{n_{\mu_i}(x)}\tilde{\Phi}_{\mu_i}(x_0)=:\mu_i^{n_{\mu_i}(x)}A_i
\end{cases}
$$
By passing to $A:=\min_iA_i$ and
$B:=\max_iB_i$ and using that
$n_{\mu_2}(x)\in[n(x)-K,n(x)+K]$ by \eqref{76}, it follows that
$$\mu_2^{n(x)}A\mu_2^{-K}\leq\tilde{\Phi}_{\mu_2}(x)\leq\mu_2^{n(x)}B\mu_2^{K}$$
and that
$$\mu_1^{n(x)}A\leq\tilde{\Phi}_{\mu_1}(x)\leq\mu_2^{n(x)}B$$
for all $x\geq x_1$. Thus there exist constants $K_1,K_2>0$ with
$$K_1\leq\frac{\tilde{\Phi}_{\mu_j}(x)}{\mu_j^{n(x)}}\leq
K_2 \text{ for all } x\geq x_1,j=1,2.$$ Now we compare
$\tilde{\Phi}_{\mu_1}$ and $\tilde{\Phi}_{\mu_2}$. We have
\begin{align*}
\tilde{\Phi}_{\mu_1}(x)&\leq K_2\mu_1^{n(x)}\\
&=K_2(\mu_2^{n(x)})^{\frac{\log(\mu_1)}{\log(\mu_2)}}\\
&\leq
K_2\frac{1}{K_1^{\frac{\log(\mu_1)}{\log(\mu_2)}}}(\tilde{\Phi}_{\mu_2}(x))^{\frac{\log(\mu_1)}{\log(\mu_2)}}\\
&=:C'\,
(\tilde{\Phi}_{\mu_2}(x))^{\frac{\log(\mu_1)}{\log(\mu_2)}}
\end{align*}
and analogously
$$\tilde{\Phi}_{\mu_1}(x)\geq
K_1\frac{1}{K_2^{\frac{\log(\mu_1)}{\log(\mu_2)}}}(\tilde{\Phi}_{\mu_2}(x))^{\frac{\log(\mu_1)}{\log(\mu_2)}}
=:c'\,
(\tilde{\Phi}_{\mu_2}(x))^{\frac{\log(\mu_1)}{\log(\mu_2)}}.$$
Let $a:=\mu_1^{\gamma_1}=\mu_2^{\gamma_2}$, i.e.
$$\gamma_i=\frac{\log a}{\log\mu_i}\text{ for }i=1,2.$$  It follows that
\begin{align*}
\tilde{\Phi}_{\mu_1}(x)^{\gamma_1}&=\tilde{\Phi}_{\mu_1}(x)^{\frac{\log
a}{\log\mu_1}}\\
&\leq C'^{\frac{\log
a}{\log\mu_1}}\tilde{\Phi}_{\mu_2}(x)^{\frac{\log
a}{\log\mu_1}\frac{\log\mu_1}{\log\mu_2}}\\
&=:C\tilde{\Phi}_{\mu_2}(x)^{\frac{\log a}{\log\mu_2}}\\
&=C\tilde{\Phi}_{\mu_2}(x)^{\gamma_2}
\end{align*}
Similarly, there exists $c>0$ such that
$$\tilde{\Phi}_{\mu_1}(x)^{\gamma_1}\geq
c\tilde{\Phi}_{\mu_2}(x)^{\gamma_2}$$ for all $x\geq x_1$. Thus,
by the definition of $\tilde{h}_{\mu_i,\gamma_i}$,
$$c\tilde{h}_{\mu_2,\gamma_2}(t)\leq\tilde{h}_{\mu_1,\gamma_1}(t)\leq
C\tilde{h}_{\mu_2,\gamma_2}(t)$$ for all $t>0$ such that
$1/t\geq x_1$.
\end{proo}
\begin{coro}\label{compare}
Let $0<\lambda_1,\lambda_2<1/e$ and $\gamma_1,\gamma_2>0$ such
that $\beta_{\lambda_1}^{\gamma_1}=\beta_{\lambda_2}^{\gamma_2}$.
Let $\Phi_{\lambda_1},\Phi_{\lambda_2}$ and
$h_{\lambda_1,\gamma_1},h_{\lambda_2,\gamma_2}$ be defined as
before. Then there exist constants $c,C>0$ such that
$$c<\frac{h_{\lambda_1,\gamma_1}(t)}{h_{\lambda_2,\gamma_2}(t)}<C$$
if $t>0$ is small enough.
\end{coro}

\subsection{Proof of the main theorem}

The proof of the main theorem is now just a combination of the results we achieved in the preceding sections. \\

\begin{proo3}
For the proof of the first statement, let $\lambda\in\Com\setminus\set{0}$. Choose $\gamma>0$ such that
\begin{equation}\label{101}
\frac{\log
g(t)}{\log\Phi_{\lambda_0}(1/t)}\geq\gamma>K_{\lambda_0}\text{ for
}t\text{ small enough}.
\end{equation}
Then $$g(t)\geq\Phi_{\lambda_0}(1/t)^\gamma$$ if
$t$ is small enough, which clearly implies that
$$\mathcal{H}^h(\mathcal{J}(E_\lambda))\geq\mathcal{H}^{h_{\lambda_0,\gamma}}(\mathcal{J}(E_\lambda)).$$
Further, \eqref{101} shows that $\beta_{\lambda_0}^\gamma>2$. Thus
$\mathcal{H}^{h_{\lambda_0,\gamma}}(\mathcal{J}(E_\lambda))=\infty$
by Theorem \ref{infmeas}, and the first statement is proven. \\
Now we prove the second statement. The hypothesis implies that there exists $\gamma_0>0$ with
\begin{equation}\label{102}
\frac{\log
g(t)}{\log\Phi_{\lambda_0}(1/t)}\leq\gamma_0<K_{\lambda_0}\text{
for }t\text{ small enough}.
\end{equation}
From \eqref{102}, it follows that
$\beta_{\lambda_0}^{\gamma_0}<2$. Now choose $\lambda\in(0,1/e)$
and $\gamma>0$ such that
$\beta_\lambda^\gamma=\beta_{\lambda_0}^{\gamma_0}$. It follows
from Corollary \ref{compare} that there exists a constant $C>0$
such that
$$\Phi_{\lambda_0}\Big(\frac{1}{t}\Big)^{\gamma_0}\leq
C\Phi_\lambda\Big(\frac{1}{t}\Big)^\gamma$$ if $t$ is small
enough. Hence
$$g(t)\leq\Phi_{\lambda_0}\Big(\frac{1}{t}\Big)^{\gamma_0}\leq
C\Phi_\lambda\Big(\frac{1}{t}\Big)^\gamma,$$ which implies that
$$\mathcal{H}^h(\mathcal{J}(E_\lambda))\leq
C\mathcal{H}^{h_{\lambda,\gamma}}(\mathcal{J}(E_\lambda))=0$$ by
Theorem \ref{zeromeas}. Thus the second statement is proven.
\end{proo3}

Note that by the remarks after the proof of Theorem \ref{infmeas} and Theorem \ref{zeromeas},
the statement remains true if we replace $\mathcal{J}(E_\lambda)$ by $I(E_\lambda)$.

\subsection{Consequences}

Finally, we show how results of Astala and Clop, as well as
Rempe, combined with Theorem \ref{uz}, can be used to
give similar estimates for Julia sets of hyperbolic and escaping
sets of arbitrary exponential maps, but - at least a priori - at the expense of the
optimal constant $K_{\lambda_0}$.
Rempe (\cite{rempe}, see also \cite{rempe2}) proved that for every $\lambda\in\Com\setminus\set{0}$,
there exists $R>0$ and a $K$-quasiconformal map
$\phi:\Com\to\Com$ such that
$$E\circ\phi=\phi\circ E_\lambda\text{ on
}I_R(E_\lambda).$$ In fact, he proved a more general
theorem, but here we only need this result. Because $\phi$ is a
conjugacy on a subset of $I(E_\lambda)$, it follows easily that
\begin{equation}\label{999}
\phi(I_R(E_\lambda))\subset I(E).
\end{equation}
Recently, Astala and Clop (oral communication) proved the following
result:
\begin{theo}\label{ascl}
Let $\gamma>0$ and $h_\gamma(t)=t^2g(1/t)^\gamma$ be a function
such that $g$ satisfies the following three assumptions:
\begin{itemize}
\item $g$ is monotonically increasing and smooth \\
\item
$\lim\limits_{t\to\infty}g(t)/t^\alpha=\begin{cases}\infty&\alpha\leq
0\\0&\alpha>0\end{cases}$
\item For each $\alpha>0$, there exists $C_\alpha>0$ and
$t_\alpha>0$ such that $$\frac{1}{C_\alpha}g(t)\leq
g(t^\alpha)\leq C_\alpha g(t)$$ whenever $t\geq t_\alpha$.
\end{itemize}
Let $\varphi:\Com\to\Com$ be a $K$-quasiconformal mapping and let
$F\subset\Com$ be compact with $\mathcal{H}^{h_\gamma}(F)=0$. Then
$$\mathcal{H}^{h_{\delta}}(\varphi(F))=0\text{ for every
}\delta<\frac{\gamma}{K}.$$
\end{theo}
It is easy to see that $\Phi$ satisfies the three
assumptions from the theorem. By Theorem \ref{zeromeas}, we have
$\mathcal{H}^{h_\gamma}(\mathcal{J}(E))=0$, in
particular
$\mathcal{H}^{h_\gamma}(\mathcal{J}(E)\cap\overline{D(0,n)})=0$
for every $n\in\mathbb{N}$. Let $\delta<\gamma/K$. Because
$\mathcal{J}(E)\cap\overline{D(0,n)}$ is a compact set and
$\phi^{-1}$ is also a $K$-quasiconformal mapping, it follows from
Theorem \ref{ascl} and $I(E)\subset\mathcal{J}(E)$ that
$$\mathcal{H}^{h_{\delta}}(\phi^{-1}(I(E)\cap\overline{D(0,n)}))=0 \text{ for every }n\in\mathbb{N}.$$
Since $$\phi^{-1}(I(E))=\bigcup_{n\in\mathbb{N}}\phi^{-1}(I(E)\cap\overline{D(0,n)}),$$
this implies $\mathcal{H}^{h_{\delta}}(\phi^{-1}(I(E)))=0.$
Now, we conclude from \eqref{999} that $\mathcal{H}^{h_{\delta}}(I_R(E_\lambda))=0,$
and Lemma \ref{irreicht} yields $\mathcal{H}^{h_{\delta}}(I(E_\lambda))=0.$
Thus we obtain
\begin{theo}\label{main2}
For every $\lambda_0\in(0,1/e)$ and
$\lambda\in\Com\setminus\set{0}$, there exist constants
$K_1,K_2>0$ (where $K_1=\log 2/\log\beta_{\lambda_0}$) with the
following property: \\ Let $h(t)=t^2g(t)$ be a gauge function.
\begin{enumerate}
\item If
$$\liminf_{t\to0}\frac{\log
g(t)}{\log\Phi_{\lambda_0}(1/t)}>K_1,$$ then
$\mathcal{H}^h(I(E_\lambda))=\infty$. The measure $\mathcal{H}^h$
is not even $\sigma$-finite on $I(E_\lambda)$.
\item If
$$\limsup_{t\to0}\frac{\log
g(t)}{\log\Phi_{\lambda_0}(1/t)}<K_2,$$ then
$\mathcal{H}^h(I(E_\lambda))=0$.
\end{enumerate}
\end{theo}
In the hyperbolic case, we get the same result for
$\mathcal{J}(E_\lambda)$ by Theorem \ref{uz}, because in this
situation, the set $\mathcal{J}(E_\lambda)\setminus I(E_\lambda)$
has zero $\mathcal{H}^{h_{\lambda_0,\gamma}}$-measure for every
$\lambda_0\in(0,1/e)$ and every $\gamma>0$. We state this result
for completeness:
\begin{theo}\label{main3}
For every $\lambda_0\in(0,1/e)$ and
$\lambda\in\Com\setminus\set{0}$ such that $E_\lambda$ is
hyperbolic, there exist constants $K_1,K_2>0$ (where $K_1=\log 2/\log\beta_{\lambda_0}$)
with the following property: \\ Let
$h(t)=t^2g(t)$ be a gauge function.
\begin{enumerate}
\item If
$$\liminf_{t\to0}\frac{\log
g(t)}{\log\Phi_{\lambda_0}(1/t)}>K_1,$$ then
$\mathcal{H}^h(\mathcal{J}(E_\lambda))=\infty$. The measure
$\mathcal{H}^h$ is not even $\sigma$-finite on
$\mathcal{J}(E_\lambda)$.
\item If
$$\limsup_{t\to0}\frac{\log
g(t)}{\log\Phi_{\lambda_0}(1/t)}<K_2,$$ then
$\mathcal{H}^h(\mathcal{J}(E_\lambda))=0$.
\end{enumerate}
\end{theo}

\end{document}